\documentclass[12pt, a4paper]{article}
\usepackage{graphicx,psfrag}
\usepackage{amsmath}
\usepackage{amssymb}
\usepackage{calligra}
\title{Efficient time stepping for the multiplicative Maxwell fluid including the Mooney-Rivlin hyperelasticity}
\author{A.V. Shutov  \\
	Lavrentyev Institute of Hydrodynamics, \\ pr. Lavrentyeva 15, 630090, Novosibirsk, Russia   \\
    Novosibirsk State University, \\ Pirogova 1, 630090, Novosibirsk, Russia   \\
	}

\date{\today}
\begin{document}
\maketitle

\begin{abstract}
A popular version of the finite strain Maxwell fluid is considered, which is based on the
multiplicative decomposition of the deformation gradient tensor. The model combines Newtonian
viscosity with hyperelasticity of Mooney-Rivlin type; it is a special case
of the viscoplasticty model proposed by Simo and Miehe (1992). A simple, efficient
and robust implicit time stepping
procedure is suggested.
Lagrangian and Eulerian versions of the algorithm are available, with equivalent properties.
 The numerical scheme is iteration free, unconditionally stable and first order accurate.
It exactly preserves the inelastic incompressibility, symmetry, positive definiteness
of the internal variables, and w-invariance.
The accuracy of the stress computations is tested using a series of numerical simulations involving
a non-proportional loading and large strain increments.
In terms of accuracy, the proposed algorithm is equivalent to the modified Euler backward method with exact
inelastic incompressibility; the proposed method is also equivalent to the classical
integration method based on exponential mapping. Since the new method
is iteration free, it is more robust and computationally efficient.
The algorithm is implemented into MSC.MARC and
a series of initial boundary value problems is solved in order to demonstrate
the usability of the numerical procedures.

\end{abstract}


\section{Introduction}

The Maxwell fluid, also known as the Maxwell material (cf. \cite{Reiner}),
is in the focus of the current research. Thanks to a special combination
of elastic and viscous (inelastic) properties, it is widespread
in the phenomenological material modelling.
Apart from obvious applications to the dynamics of viscous fluids \cite{BalanTsakmakis},
different versions of the finite strain Maxwell model are employed
for solid materials as well. A parallel connection of the Maxwell fluid
with idealized elastic elements allows one to model the behaviour of
polymers \cite{LionAM, ReeseHab, Hasanpour, Kleuter, Diani, Ghobadi, Johlitz}, biological
tissues \cite{Holzapfel, GassFor, Koolstra, Vassoler}, explosives \cite{Buechler} and
other types of viscoelastic materials.
Plasticity models
with a nonlinear kinematic hardening may operate with backstresses;
the evolution of the backstresses can be described using the Maxwell element as well
\cite{LionIJP, DettRes, Vladimirov, Feigen, ShutovKuprin}.
A generalization of this approach to plasticity with the yield surface distortion
is reported in \cite{ShutovPaKr, ShutovIhle}.
Some advanced models of
shape memory alloys \cite{Helm1, Helm3} and anisotropic creep \cite{Shutov2017}
include the Maxwell fluid as an important constituent.
Interestingly, minor modifications of the Maxwell fluid are also
used in production of mechanics-based
computer animations \cite{Ram}.
Applications of the Maxwell fluid to the analysis of
geological structures are reported in \cite{SimoMeschke, PericCr}.
Due to the high prevalence of the Maxwell concept in the phenomenological material modeling,
specialized numerical algorithms with an enhanced robustness and efficiency are becoming
increasingly important. High computational efficiency of the utilized algorithms
may become especially important
for real-time simulations, like FEM simulations
of surgical operations \cite{Misra}.

Different mathematical formulations of
the finite strain Maxwell fluid are currently in use (cf. \cite{ShutovSeven}). We advocate here the
multiplicative approach of Simo and Miehe \cite{SimMieh},
due to its numerous advantages like
thermodynamic consistency, exact hyperelasticity, absence of spurious shear oscillations,
pure isochoric-volumetric split \cite{ShutovSeven}, w-invariance under
the isochoric change of the reference configuration \cite{WInvariance} and
exponential stability of the solution with respect to the perturbation of the initial data
\cite{ShutovStability}. Within the approach of Simo and Miehe, the Newtonian viscosity can be combined with
different types of isotropic elastic potentials in a thermodynamically consistent manner.

In computational practice, the implemented time step size can be larger than the relaxation time of a Maxwell element.
For that reason, implicit integration of the six-dimensional evolution equation is needed.
After the discretization in time, a system of six nonlinear equations with respect to six unknowns is obtained.
Using the spectral decomposition, the number of unknowns can be reduced to three \cite{Simo, ReeseG, HolLoug, Haider}.
Typically, these nonlinear systems of equations are solved iteratively using the Newton-Raphson method
or its modifications \cite{ReeseG, HartmannHabil, PericCr, Helm2, Kleuter, ShutovKrKoo, Vladimirov, Hasanpour,
HolLoug, Rauchs, Lejeunes, Vassoler}. Obviously, such an iterative procedure is less robust and less efficient than a procedure
based on the closed form solution of the discretized system of equations.

In the important case of the neo-Hookean potential, an explicit update formula for the implicit time
stepping was proposed in \cite{ShutovLandgraf}. In \cite{ShutovLandgraf} it was used to model finite strain viscoelastic behaviour
of a rubber-like material. Further extension of this algorithm to the case of the Yeoh hyperelasticity
was presented in \cite{LandgrafShutov}.  In \cite{Johlitz, Ghobadi} the algorithm was extended to the thermo-mechanical case.
In \cite{Silbermann, ShutovDamage} the algorithm was used to stabilize numerical simulations
within a hybrid explicit/implicit procedure for finite strain plasticity with a nonlinear kinematic hardening.
In \cite{Schueler} the algorithm was used to model the viscoelastic behavior of a bituminous binding agent.
The closed form solution reported in \cite{ShutovLandgraf} is implemented as a part
of efficient implicit procedures for finite strain viscoplasticity \cite{Shutov2016} and
finite strain creep \cite{Shutov2017}.

In the current study, a new explicit update formula is suggested for a more
general case of the Mooney-Rivlin hyperelasticity.
The explicit solution is obtained by exploiting the properties
of the underlying constitutive equations.
The previously reported update formula for the
neo-Hookean potential is covered by this solution as a special case. The new algorithm
is unconditionally stable and first order accurate.\footnote{For the discussion of
higher order methods, the reader is referred to \cite{ReeseHab, Hartmann2010, Eidel, EidelStumpf}.}
It exactly preserves
the inelastic incompressibility, symmetry and positive definiteness of the
internal tensor-valued variable; it also preserves the
w-invariance under the change of the reference configuration. A slight modification
of the method is suggested to enforce the symmetry of the consistent tangent operator.
Concerning the accuracy of the stress computation, the new algorithm is equivalent to the
modified Euler backward method (MEBM) and exponential method (EM).
The symmetry of the consistent tangent operator is tested in a series of computations.
In order to demonstrate the applicability of the numerical procedure, an initial
boundary value problem is solved in MSC.MARC.

We conclude the introduction by setting up the notation.
A coordinate free tensor formalism is used here.
Bold-faced symbols denote first- and second-rank tensors in $\mathbb{R}^3$. For instance,
$\mathbf{1}$ stands for the second-rank identity tensor.
The deviatoric part of a tensor is denoted by
$\mathbf A^{\text{D}} := \mathbf A - \frac{1}{3} \text{tr}(\mathbf A) \mathbf 1$, where
$\text{tr}(\mathbf A)$ is the trace.
The overline $\overline{(\cdot)}$ denotes the unimodular
part of a tensor:
\begin{equation}\label{BarDef}
\overline{\mathbf{A}}:=(\det \mathbf{A})^{-1/3} \mathbf{A}.
\end{equation}

\section{Finite strain Maxwell fluid according to Simo and Miehe}

The considered Maxwell model is a special case of the finite strain viscoplasticity theory proposed by Simo and Miehe \cite{SimMieh};
it has the same structure as the well known model of associative elastoplasticity introduced by Simo in \cite{Simo}.
A referential (Lagrangian) formulation of this model was considered later by Lion in \cite{LionAM}.
The spatial (Eulerian) constitutive equations proposed by
Simo \& Miehe were used by Reese and Govindjee in the fundamental study \cite{ReeseG}
and by many others (see, for instance, \cite{HubTsak, Nedjar, PericCr,
Kleuter, Hasanpour, HolLoug, Rauchs, Lejeunes}).

\subsection{Formulation on the reference configuration}

Here we follow the presentation of Lion \cite{LionAM}. Let $\mathbf F$ be the deformation gradient which maps
the local reference configuration  $\tilde{\mathcal{K}}$
to the current configuration $\mathcal{K}$.
We consider the multiplicative decomposition of the deformation gradient $\mathbf F$
into the elastic part $\hat{\mathbf F}_{\text{e}}$ and the inelastic part $\mathbf F_{\text{i}}$
\footnote{In the context of viscoelasticity, this multiplicative split is known
as the Sidoroff decomposition \cite{Sidoroff}.}
\begin{equation*}\label{split1}
\mathbf F = \hat{\mathbf F}_{\text{e}} \mathbf F_{\text{i}}.
\end{equation*}
This decomposition implies the so-called stress-free intermediate configuration
$\hat{\mathcal{K}}$.
Along with the classical right Cauchy-Green tensor
${\mathbf C}: =\mathbf F^{\text{T}} \mathbf F$ we consider
the inelastic right Cauchy-Green tensor ${\mathbf C}_{\text{i}}$ and
the elastic right Cauchy-Green tensor $\hat{\mathbf C}_{\text{e}}$
\begin{equation}\label{intIntv}
{\mathbf C}_{\text{i}} :=\mathbf F^{\text{T}}_{\text{i}}
\mathbf F_{\text{i}}, \quad \hat{\mathbf C}_{\text{e}}  :=
\hat{\mathbf F}_{\text{e}}^{\text{T}} \hat{\mathbf F}_{\text{e}}.
\end{equation}
The inelastic velocity gradient $\hat{\mathbf L}_{\text{i}} $ and the corresponding
covariant Oldroyd derivative are defined as follows
\begin{equation*}\label{altol}
\hat{\mathbf L}_{\text{i}} := \dot{\mathbf F}_{\text{i}} \mathbf F^{-1}_{\text{i}}, \quad
\stackrel{\triangle} {(\cdot)} := \frac{d}{d t}(\cdot) +
\hat{\mathbf L}^{\text{T}}_{\text{i}}
(\cdot) + (\cdot) \hat{\mathbf L}_{\text{i}}.
\end{equation*}
Here, the superimposed dot denotes the material time rate.
The inelastic Almansi strain tensor $\hat{\mathbf \Gamma}_{\text{i}} $ and the inelastic
strain rate tensor $\hat{\mathbf D}_{\text{i}} $ are introduced through
 \begin{equation*}\label{defgae}
\hat{\mathbf \Gamma}_{\text{i}} := \frac{\displaystyle 1}{\displaystyle 2}
(\mathbf 1 - \mathbf F_{\text{i}}^{-\text{T}} \mathbf F_{\text{i}}^{-1}), \quad
\hat{\mathbf D}_{\text{i}} := \frac{1}{2}
(\hat{\mathbf L}_{\text{i}} + \hat{\mathbf L}^{\text{T}}_{\text{i}} ).
\end{equation*}
After some straightforward computations, we arrive at
\begin{equation}\label{DOldrDer}
\hat{\mathbf D}_{\text{i}} = \stackrel{\triangle}{\hat{\mathbf \Gamma}}_{\text{i}}.
\end{equation}

By $\mathbf T$ denote the Cauchy stress tensor (true stresses).
The Kirchhoff stress $\mathbf S$,
the 2nd Piola-Kirchhoff stress $\hat{\mathbf S}$ operating on the stress-free configuration $\hat{\mathcal{K}}$,
and the classical 2nd Piola-Kirchhoff stress $\tilde{\mathbf T}$ operating on the
reference configuration  $\tilde{\mathcal{K}}$ are now introduced:
\begin{equation}\label{Kirch}
\mathbf S  : = (\text{det} \mathbf F) \mathbf T, \quad
\hat{\mathbf S} : = {\hat{\mathbf F}_{\text{e}}^{-1}} \mathbf S {\hat{\mathbf F}_{\text{e}}^{-\text{T}}}, \quad
\tilde{\mathbf T} : = {\mathbf F^{-1}} \mathbf S {\mathbf F^{-\text{T}}}.
\end{equation}
Let $\psi$ be the Helmholz free energy per unit mass. In this work it is given by the Mooney-Rivlin potential
\begin{equation}\label{spec1}
\rho_{\scriptscriptstyle \text{R}}  \psi(\hat{\mathbf{C}}_{\text{e}})=
\frac{c_{10}}{2} \big( \text{tr} \overline{\hat{\mathbf{C}}_{\text{e}}} - 3 \big) +
\frac{c_{01}}{2} \big( \text{tr} \overline{\hat{\mathbf{C}}^{-1}_{\text{e}}} - 3 \big),
\end{equation}
where $\rho_{\scriptscriptstyle \text{R}}$ is the mass density in the reference
configuration, $c_{10}$ and $c_{01}$ are the shear moduli; the overline $\overline{(\cdot)}$
denotes the unimodular part of a tensor (recall \eqref{BarDef}).\footnote{The
neo-Hookean potential is obtained when $c_{01}=0$; this special
 case was already considered in \cite{ShutovLandgraf}.}
A hyperelastic stress-strain relation is considered
on the stress-free configuration:
\begin{equation}\label{potent}
\hat{\mathbf S} = 2 \rho_{\scriptscriptstyle \text{R}}
\frac{\displaystyle \partial \psi(\hat{\mathbf C}_{\text{e}})}
{\displaystyle \partial \hat{\mathbf C}_{\text{e}}} \quad \Rightarrow  \quad
\hat{\mathbf S} = \Big(c_{10} \ \hat{\mathbf C}_{\text{e}} -
c_{01} \ \hat{\mathbf C}^{-1}_{\text{e}} \Big)^{\text{D}} \hat{\mathbf C}^{-1}_{\text{e}}.
\end{equation}

The Clausius-Duhem inequality requires that
the specific internal dissipation $\delta_{\text{i}}$ remains non-negative.
For simplicity, we consider isothermal processes here and the
Clausius-Duhem inequality takes the reduced form, also known as the Clausius-Planck inequality
\begin{equation}\label{cld}
\delta_{\text{i}} := \frac{1}{\rho_{\scriptscriptstyle \text{R}}} \tilde{\mathbf T} : \dot{\mathbf E} - \dot{\psi} \geq 0,
\end{equation}
where $\mathbf E : = \frac{\displaystyle 1}{\displaystyle 2} (\mathbf C - \mathbf 1)$ stands for the
Green strain tensor.
Using $\eqref{potent}_1$ and taking the isotropy of the free energy function into account,
this inequality is reduced to
\begin{equation}\label{cld2}
\rho_{\scriptscriptstyle \text{R}} \delta_{\text{i}} =
(\hat{\mathbf{C}}_{\text{e}}  \hat{\mathbf S}) : \stackrel{\triangle} {\hat{\mathbf{\Gamma}}}_{\text{i}}  \geq 0.
\end{equation}
The following flow rule is postulated so that
\eqref{cld2} holds for arbitrary
mechanical loadings (cf. \cite{LionAM})
\begin{equation}\label{evol}
\stackrel{\triangle}{\hat{\mathbf{\Gamma}}}_{\text{i}} = \frac{1}{2 \eta}
(\hat{\mathbf{C}}_{\text{e}}  \hat{\mathbf S})^{\text{D}},
\end{equation}
$\eta > 0$ is a material parameter (Newtonian viscosity).
In view of \eqref{DOldrDer}, an equivalent formulation of this flow
rule is given by
\begin{equation}\label{evol2}
\hat{\mathbf{D}}_{\text{i}} = \frac{1}{2 \eta}
(\hat{\mathbf{C}}_{\text{e}}  \hat{\mathbf S})^{\text{D}}.
\end{equation}
Both \eqref{evol} and \eqref{evol2} imply that $\text{tr} (\stackrel{\triangle}{\hat{\mathbf{\Gamma}}}_{\text{i}}) =
\text{tr} (\hat{\mathbf{D}}_{\text{i}}) =0$.
The inelastic flow is thus incompressible: $\det \mathbf{F}_{\text{i}} = const$.

Let us transform the constitutive equations to the reference configuration.
First, the free energy \eqref{spec1} is represented as a function of $\mathbf C {\mathbf C_{\text{i}}}^{-1}$:
\begin{equation}\label{freeen2}
\psi=\psi(\mathbf C {\mathbf C_{\text{i}}}^{-1})=
\frac{c_{10}}{2 \rho_{\scriptscriptstyle \text{R}} } \big(
\text{tr} \overline{\mathbf C {\mathbf C_{\text{i}}}^{-1}} - 3 \big) +
\frac{c_{01}}{2 \rho_{\scriptscriptstyle \text{R}} } \big(
\text{tr} \overline{ \mathbf C_{\text{i}} \mathbf C^{-1}} - 3 \big).
\end{equation}
Using \eqref{Kirch} and $\eqref{potent}_1$, one obtains for the 2nd Piola-Kirchhoff
stress tensor
\begin{equation}\label{UseinAppD}
\tilde{\mathbf T}  =
2 \rho_{\scriptscriptstyle \text{R}}
\frac{\displaystyle \partial \psi(\mathbf C {\mathbf C_{\text{i}}}^{-1})}
{\displaystyle \partial \mathbf{C}}\big|_{\mathbf C_{\text{i}} =
\text{const}}.
\end{equation}
Substituting \eqref{freeen2} into this, we arrive at
\begin{equation}\label{2PKd}
\tilde{\mathbf T}  = \mathbf C^{-1} (c_{10} \overline{\mathbf C} \mathbf C_{\text{i}}^{-1}
- c_{01} \mathbf C_{\text{i}} \overline{\mathbf C^{-1}} )^{\text{D}}.
\end{equation}
Next, due to the isotropy of the free energy function we have
\begin{equation}\label{trace}
\text{tr} (\hat{\mathbf{C}}_{\text{e}}  \hat{\mathbf S}) =
\text{tr} (\mathbf C \tilde{\mathbf T}).
\end{equation}
Combining this with \eqref{Kirch} we obtain
\begin{equation}\label{puba2}
{\mathbf F}_{\text{i}}^{\text{T}} (\hat{\mathbf{C}}_{\text{e}}  \hat{\mathbf S})^{\text{D}} {\mathbf F}_{\text{i}} =
\mathbf C \tilde{\mathbf T} {\mathbf C}_{\text{i}}  - \frac{1}{3} \text{tr} (\hat{\mathbf{C}}_{\text{e}}  \hat{\mathbf S})
\ {\mathbf C}_{\text{i}} \stackrel{\eqref{trace}}{=}
\big( \mathbf C \tilde{\mathbf T} \big)^{\text{D}} \mathbf C_{\text{i}}.
\end{equation}
Applying a pull-back transformation to the evolution equation \eqref{evol} we obtain
\begin{equation}\label{puba4}
\dot{\mathbf C}_{\text{i}} = 2 {\mathbf F}_{\text{i}}^{\text{T}}
\stackrel{\triangle}{\hat{\mathbf{\Gamma}}}_{\text{i}} {\mathbf F}_{\text{i}} \stackrel{\eqref{evol}}{=}
 \frac{1}{\eta} {\mathbf F}_{\text{i}}^{\text{T}} (\hat{\mathbf{C}}_{\text{e}}  \hat{\mathbf S})^{\text{D}} {\mathbf F}_{\text{i}}
  \stackrel{\eqref{puba2}}{=}  \frac{1}{\eta} \big( \mathbf C \tilde{\mathbf T} \big)^{\text{D}} \mathbf C_{\text{i}}.
\end{equation}
This flow rule is valid for arbitrary isotropic
free energy functions.
Substituting \eqref{2PKd} into \eqref{puba4}, we obtain the evolution
equation pertaining to the Mooney-Rivlin potential
\begin{equation}\label{puba44}
\dot{\mathbf C}_{\text{i}} = \frac{1}{\eta}  \big( c_{10} \overline{\mathbf C}
\mathbf C^{-1}_{\text{i}} - c_{01}
\mathbf C_{\text{i}}  \overline{\mathbf C^{-1}} \big)^{\text{D}} \mathbf C_{\text{i}}.
\end{equation}
Finally, the system of constitutive equations \eqref{2PKd} and \eqref{puba44} is closed by
specifying initial conditions
\begin{equation}\label{initCond}
\mathbf C_{\text{i}}|_{t=t^0} = \mathbf C_{\text{i}}^0.
\end{equation}

The exact solution of \eqref{puba44} exhibits
the following geometric property
\begin{equation}\label{geopro}
\mathbf{C}_{\text{i}}(t) \in \mathbb{M} \quad \text{if} \quad
\mathbf{C}^0_{\text{i}} \in \mathbb{M},
\end{equation}
where the manifold $\mathbb{M}$ is a set of symmetric unimodular tensors
\begin{equation*}\label{UnimodulM}
\mathbb{M} := \big\{ \mathbf A \in Sym: \text{det} \mathbf A =1 \big\}.
\end{equation*}
In other words, we a dealing
with a differential equation of the manifold.
A substantial consequence of \eqref{geopro} is that $\mathbf{C}_{\text{i}}$
remains positive definite.
The positive definiteness is important since  $\mathbf{C}_{\text{i}}$ represents
a certain metric in $\mathbb{R}^3$.

\subsection{Spatial formulation}

Similar to \cite{ShutovLandgraf}, we show that the model presented in this section
is indeed the model of Simo and Miehe \cite{SimMieh}, which was originally formulated
in the spatial (Eulerian) description.
First,  for the inelastic strain rate we have a purely kinematic relation
\begin{equation}\label{SimMie}
2 {\hat{\mathbf{D}}}_{\text{i}} = {\mathbf F}_{\text{i}}^{-\text{T}}
\dot{\mathbf C}_{\text{i}} {\mathbf F}_{\text{i}}^{-1}.
\end{equation}
Since the elastic
potential $\psi(\hat{\mathbf{C}}_{\text{e}})$ is isotropic,
the tensors $\hat{\mathbf S}$ and $\hat{\mathbf{C}}_{\text{e}}$ are coaxial.
Thus, the Mandel tensor $ \hat{\mathbf{C}}_{\text{e}} \hat{\mathbf S}$ and the tensor $\hat{\mathbf{C}}_{\text{e}}$ are coaxial as well.
The evolution equation \eqref{evol2} implies that ${\hat{\mathbf{D}}}_{\text{i}}$ is coaxial with $\hat{\mathbf{C}}_{\text{e}}$.
Therefore, these tensors commute: ${\hat{\mathbf{D}}}_{\text{i}} \hat{\mathbf{C}}_{\text{e}} = \hat{\mathbf{C}}_{\text{e}} {\hat{\mathbf{D}}}_{\text{i}}$
and we have
\begin{equation*}\label{SimMie3456643}
{\hat{\mathbf{D}}}_{\text{i}} = \hat{\mathbf{C}}_{\text{e}} {\hat{\mathbf{D}}}_{\text{i}} \hat{\mathbf{C}}^{-1}_{\text{e}}.
\end{equation*}
Taking this into account, we rewrite \eqref{evol2} as follows
\begin{equation}\label{SimMie2}
2 \hat{\mathbf{C}}_{\text{e}} \hat{\mathbf{D}}_{\text{i}}
\hat{\mathbf{C}}_{\text{e}}^{-1} = \frac{1}{\eta}
(\hat{\mathbf{C}}_{\text{e}}  \hat{\mathbf S})^{\text{D}}.
\end{equation}
Substituting \eqref{SimMie} into \eqref{SimMie2} and taking into account
that $\frac{d}{d t} ({\mathbf C}_{\text{i}}^{-1} ) = -
{\mathbf C}_{\text{i}}^{-1} \dot{{\mathbf C}}_{\text{i}} {\mathbf C}_{\text{i}}^{-1}$, we obtain
\begin{equation}\label{SimMie3}
-\hat{\mathbf{F}}^{\text{T}}_{\text{e}} \mathbf{F}
\frac{d}{d t} ({\mathbf C}_{\text{i}}^{-1} ) \mathbf{F}^{\text{T}}
\hat{\mathbf{F}}_{\text{e}}^{-\text{T}}
\hat{\mathbf{F}}_{\text{e}}^{-1} \hat{\mathbf{F}}_{\text{e}}^{-\text{T}}   = \frac{1}{\eta}
(\hat{\mathbf{C}}_{\text{e}}  \hat{\mathbf S})^{\text{D}}.
\end{equation}
For the Eulerian description it is convenient to introduce the elastic left Cauchy-Green tensor
$\mathbf{B}_{\text{e}} := \hat{\mathbf{F}}_{\text{e}} \hat{\mathbf{F}}_{\text{e}}^{\text{T}}$. Next we
note that the Kirchhoff stress is related to the Mandel tensor through the similarity relation
\begin{equation}\label{SimMie4}
\mathbf{S}^{\text{D}} = \hat{\mathbf{F}}_{\text{e}}^{-\text{T}}
(\hat{\mathbf{C}}_{\text{e}}  \hat{\mathbf S})^{\text{D}}\hat{\mathbf{F}}_{\text{e}}^{\text{T}}.
\end{equation}
Pre-multiplying \eqref{SimMie3} with $\hat{\mathbf{F}}^{-\text{T}}_{\text{e}}$, post-multiplying it
with $\hat{\mathbf{F}}^{\text{T}}_{\text{e}}$ and taking \eqref{SimMie4} into account we have
\begin{equation}\label{SimMie5}
-\mathbf{F}
\frac{d}{d t} ({\mathbf C}_{\text{i}}^{-1} ) \mathbf{F}^{\text{T}}
\mathbf{B}_{\text{e}}^{-1}
 = \frac{1}{\eta} \mathbf S^{\text{D}}.
\end{equation}
Introducing the covariant Oldroyd rate (which is effectively a Lie derivative)
\begin{equation}\label{OldrCovar}
\mathfrak{O} (\mathbf{A}) =
\text{\calligra{L}}_{v}  (\mathbf{A}) : =
\mathbf{F} \frac{d}{dt} (\mathbf{F}^{-1} \mathbf{A} \mathbf{F}^{-\text{T}})  \mathbf{F}^{\text{T}}
= \dot{\mathbf{A}} -  \mathbf{L}  \mathbf{A} -  \mathbf{A}  \mathbf{L}^{\text{T}},
\end{equation}
we obtain the following kinematic relation
\begin{equation*}\label{OldrCovar2}
\mathfrak{O} (\mathbf{B}_{\text{e}}) =
\text{\calligra{L}}_{v}  (\mathbf{B}_{\text{e}}) = \mathbf{F}
\frac{d}{d t} ({\mathbf C}_{\text{i}}^{-1} ) \mathbf{F}^{\text{T}}.
\end{equation*}
Using it, the flow rule \eqref{SimMie5} takes the well-known form
\begin{equation}\label{SimMie82}
-\text{\calligra{L}}_{v}  (\mathbf{B}_{\text{e}})
\mathbf{B}_{\text{e}}^{-1}
 = \frac{1}{\eta} \mathbf S^{\text{D}}.
\end{equation}
This equation coincides with the flow rule considered by Simo and Miehe \cite{SimMieh}
(see equations (2.19a) and (2.26) in \cite{SimMieh}); it was also implemented by
Reese and Govindjee in \cite{ReeseG}. The Kirchhoff stress can be computed through
\begin{equation}\label{Murnag0}
\mathbf{S} =  2 \rho_{\scriptscriptstyle \text{R}}
\frac{\displaystyle \partial \psi(\mathbf{B}_\text{e})}
{\displaystyle \partial \mathbf{B}_\text{e}} \mathbf{B}_\text{e}.
\end{equation}
In the case of the Mooney-Rivlin strain energy \eqref{freeen2}, we have
$\mathbf S = \mathbf S^{\text{D}} = c_{10} (\overline{\mathbf{B}_{\text{e}}})^{\text{D}}
-c_{01} (\overline{\mathbf{B}^{-1}_{\text{e}}})^{\text{D}} $;
the evolution equation \eqref{SimMie82} is then specified to
\begin{equation}\label{SimMie8}
-\text{\calligra{L}}_{v}  (\mathbf{B}_{\text{e}})
\mathbf{B}_{\text{e}}^{-1}
 = \frac{1}{\eta} (c_{10} \overline{\mathbf{B}_{\text{e}}} - c_{01} \overline{\mathbf{B}^{-1}_{\text{e}}})^{\text{D}}.
\end{equation}

\emph{Remark 1} \\
The flow rule \eqref{SimMie82} can be obtained by other arguments as well
(cf. the derivation of Eq. (81) in reference \cite{LatorreMontans}).
Moreover, an alternative spatial formulation of the model
\eqref{SimMie82}-\eqref{Murnag0} can be derived from the additive decomposition
of the strain rate tensor (see Appendix A).
The abundance of different formulations of this model proposed by different authors
indicates the importance of this particular model.

\section{Time stepping algorithm}

\subsection{Explicit update formula in the Lagrangian formulation}

Consider a generic time interval $(t_n, t_{n+1})$ with the time step size $\Delta t:= t_{n+1} - t_n > 0$.
By ${}^n \mathbf{C}_{\text{i}}$ and ${}^{n+1} \mathbf{C}_{\text{i}}$ denote numerical
solutions at $t_n$ and $t_{n+1}$, respectively.
Assume that the current deformation gradient ${}^{n+1} \mathbf F$
and the previous inelastic Cauchy-Green tensor ${}^n \mathbf{C}_{\text{i}} \in \mathbb{M}$ are given.
The unknown ${}^{n+1} \mathbf{C}_{\text{i}} \in \mathbb{M}$
is computed by implicit integration of the evolution equation \eqref{puba44}.
First, we consider the classical Euler backward discretization of \eqref{puba44}:
\begin{equation}\label{ClassicalEBM}
{}^{n+1}\mathbf C^{(\text{EBM})}_{\text{i}} = {}^{n}\mathbf C_{\text{i}} + \frac{\Delta t}{\eta}  \Big( c_{10} \overline{{}^{n+1} \mathbf C} \
({}^{n+1} \mathbf C^{(\text{EBM})}_{\text{i}})^{-1} - c_{01}
{}^{n+1} \mathbf C^{(\text{EBM})}_{\text{i}} \  \overline{{}^{n+1} \mathbf C^{-1}} \Big)^{\text{D}} \ {}^{n+1} \mathbf C_{\text{i}}.
\end{equation}
Unfortunately, due to its linear structure, the Euler backward method violates the incompressibility
restriction which leads to error accumulation \cite{ShutovStability}. In order to enforce the incompressibility, a correction term
$\tilde{\varphi} \ {}^{n+1} \mathbf C_{\text{i}}$ is added to the
right-hand side of \eqref{ClassicalEBM} (cf. \cite{Shutov2016})
\begin{equation}\label{ModEBM}
{}^{n+1}\mathbf C_{\text{i}} = {}^{n}\mathbf C_{\text{i}} + \frac{\Delta t}{\eta}  \big( c_{10} \overline{{}^{n+1} \mathbf C} \
({}^{n+1} \mathbf C_{\text{i}})^{-1} - c_{01}
{}^{n+1} \mathbf C_{\text{i}} \  \overline{{}^{n+1} \mathbf C^{-1}} \big)^{\text{D}} \
{}^{n+1} \mathbf C_{\text{i}} + \tilde{\varphi} \ {}^{n+1} \mathbf C_{\text{i}},
\end{equation}
where the unknown correction $\tilde{\varphi} \in \mathbb{R}$ should be
 defined from the additional equation $\det({}^{n+1}\mathbf C_{\text{i}}) = 1$.
Next, using the definition of the deviatoric part, \eqref{ModEBM} takes the equivalent form
\begin{equation}\label{ModEBM2}
\varphi \ {}^{n+1}\mathbf C_{\text{i}} = {}^{n}\mathbf C_{\text{i}} +
\frac{\Delta t}{\eta}   c_{10} \overline{{}^{n+1} \mathbf C} \
 - c_{01} \frac{\Delta t}{\eta}
{}^{n+1} \mathbf C_{\text{i}} \  \overline{{}^{n+1} \mathbf C^{-1}} \
{}^{n+1} \mathbf C_{\text{i}},
\end{equation}
where $\varphi \in \mathbb{R}$ should be defined from the incompressibility condition.
For brevity, introduce notation
\begin{equation}\label{ModEBM4}
\mathbf X : = \overline{{}^{n+1} \mathbf C^{-1/2}} \ {}^{n+1}\mathbf C_{\text{i}} \
\overline{{}^{n+1} \mathbf C^{-1/2}}, \quad
\mathbf A : = \overline{{}^{n+1} \mathbf C^{-1/2}} \big(
 {}^{n}\mathbf C_{\text{i}} + \frac{\Delta t}{\eta}   c_{10} \overline{{}^{n+1} \mathbf C} ) \overline{{}^{n+1} \mathbf C^{-1/2}}, \quad
 \varepsilon := c_{01} \frac{\Delta t}{\eta}.
\end{equation}
Multiplying both sides of \eqref{ModEBM2} with $\overline{{}^{n+1} \mathbf C^{-1/2}}$ on the left
and right, we arrive at quadratic equation with respect to $\mathbf X$
\begin{equation}\label{ModEBM3}
\varphi \ \mathbf X = \mathbf{A}
 - \varepsilon \ \mathbf X^{2}.
\end{equation}
Here, $\mathbf{A}$ and $\varepsilon$ are known; the unknown $\varphi$ is defined using the incompressibility condition
\begin{equation}\label{ModEBM5}
\det(\mathbf{X}) = 1.
\end{equation}
Recall that, for mechanical reasons, ${}^{n+1} \mathbf C_{\text{i}}$ must be positive definite.
Therefore, the physically reasonable $\mathbf X$ is positive definite as well.
Thus, the correct solution of \eqref{ModEBM3} is given by
\begin{equation}\label{ReducedEvolut5}
\mathbf X = \frac{\displaystyle
 1}{\displaystyle 2 \varepsilon} \Big[ -\varphi \textbf{1} + \big(\varphi^2 \textbf{1} +  4 \varepsilon
  \textbf{A} \big)^{1/2} \Big].
\end{equation}
Unfortunately, due to the round-off errors this relation yields unsatisfactory
results if evaluated step by step.\footnote{
Interestingly, a similar problem appeared
in the context of finite strain plasticity with neo-Hookean potentials, discussed in \cite{Shutov2016}.}
Indeed, the relation in the square bracket on the right-hand side of \eqref{ReducedEvolut5}
can be computed accurately up to a machine precision. In the case of small $\varepsilon$, however,
the corresponding error is then multiplied with the big number $1/\varepsilon$.
In order to resolve this issue, \eqref{ReducedEvolut5}
needs to be transformed.
Using the well-known identity $(\mathbf X^{1/2} -\mathbf 1) (\mathbf X^{1/2} + \mathbf 1) = \mathbf X - \mathbf 1$,
\eqref{ReducedEvolut5} can be re-written in the form
\begin{equation}\label{ReducedEvolut5notCont}
\mathbf X = 2 \textbf{A} \Big[  \Big( \varphi^2 \textbf{1} +
4 \varepsilon \textbf{A} \Big)^{1/2} + \varphi \textbf{1} \Big]^{-1}.
\end{equation}
This relation is computationally advantageous since it is free from the
previously described error magnification.

Now we need to compute the unknown $\varphi$, which
should be estimated using the incompressibility condition $\det (\mathbf X) =1$.
As shown in Appendix B, a simple formula can be obtained using the
perturbation method for small $\varepsilon$:
\begin{equation}\label{EstimOfz}
\varphi = \varphi_0 - \frac{\text{tr} \mathbf A}{3 \varphi_0} \varepsilon + O ( \varepsilon^2),
\quad \text{where} \
\varphi_0 :=  \big(\det \mathbf A \big)^{1/3}.
\end{equation}
Neglecting the terms $O ( \varepsilon^2)$, a reliable procedure is obtained,
which is exact for $c_{01}=0$. Moreover, as will be seen in the next sections, this procedure is accurate
and robust even for finite values of $\varepsilon$.
After $\varphi$ is computed, $\mathbf{X}$ is evaluated through \eqref{ReducedEvolut5notCont}.
Further, $\mathbf C_{\text{i}} $ is updated:
\begin{equation}\label{UpdateCi}
{}^{n+1} \mathbf C^*_{\text{i}} := \overline{{}^{n+1}
\mathbf C^{1/2}}  \ \mathbf X \ \overline{{}^{n+1} \mathbf C^{1/2}}.
\end{equation}
Since the variable $\varphi$ is not computed exactly, the
incompressibility condition can be violated.
To enforce it, a final correction step is needed:
\begin{equation}\label{UpdateCi_correct}
{}^{n+1} \mathbf C_{\text{i}} = \overline{{}^{n+1} \mathbf C^*_{\text{i}}} =
\overline{ {}^{n+1} \mathbf C^{1/2}  \ \mathbf X \ {}^{n+1} \mathbf C^{1/2}}.
\end{equation}
Finally, the stress is updated using \eqref{2PKd}. We call this procedure
\emph{iteration free Euler backward method} (IFEBM).
It is summarized in Table \ref{tab0}.

\begin{table}[h!]
\caption{IFEBM on the reference configuration}
\begin{center}
\begin{tabular}{|l |}
\hline
input: ${}^{n+1} \mathbf C$, ${}^{n} \mathbf C_{\text{i}}$  \\
output: ${}^{n+1} \mathbf C_{\text{i}}$, ${}^{n+1} \tilde{\mathbf T}$  \\
1: $\mathbf A  = \overline{{}^{n+1} \mathbf C^{-1/2}} \big(
 {}^{n}\mathbf C_{\text{i}} + \frac{\Delta t}{\eta}   c_{10}
 \overline{{}^{n+1} \mathbf C} ) \overline{{}^{n+1} \mathbf C^{-1/2}}$  \\
2: $\varepsilon = c_{01} \frac{\Delta t}{\eta}$    \\
3: $\varphi_0 =  \big(\det \mathbf A \big)^{1/3}$ \\
4: $\varphi = \varphi_0 - \frac{\text{tr} \mathbf A}{3 \varphi_0} \varepsilon$  \\
5: $\mathbf X = 2 \textbf{A} \Big[  \Big( \varphi^2 \textbf{1} +
4 \varepsilon \textbf{A} \Big)^{1/2} + \varphi \textbf{1} \Big]^{-1}$  \\
6: ${}^{n+1} \mathbf C_{\text{i}}  =
\overline{ {}^{n+1} \mathbf C^{1/2}  \ \mathbf X \ {}^{n+1} \mathbf C^{1/2}}$  \\
7: ${}^{n+1} \tilde{\mathbf T}  = {}^{n+1} \mathbf C^{-1}
(c_{10} \overline{{}^{n+1} \mathbf C} {}^{n+1} \mathbf C_{\text{i}}^{-1}
- c_{01} {}^{n+1} \mathbf C_{\text{i}} \overline{{}^{n+1} \mathbf C^{-1}} )^{\text{D}}$    \\  \hline
\end{tabular} \\
\end{center}
\label{tab0}
\end{table}

For the accuracy analysis, along with the IFEBM we consider its modification.
Toward that end we introduce the following functions:
\begin{equation}\label{Resid}
\mathbf{X}(\varphi):=
2 \textbf{A} \Big[  \Big( \varphi^2 \textbf{1} +
4 \varepsilon \textbf{A} \Big)^{1/2} + \varphi \textbf{1} \Big]^{-1}, \quad
R(\varphi) : = \det(\mathbf{X}(\varphi)) - 1.
\end{equation}
Next, we consider two Newton iterations for the equation $R(\varphi) = 0$ with the initial
approximation given by \eqref{EstimOfz}
\begin{equation}\label{IterProcess}
\varphi^{(0)} : = \varphi_0 - \frac{\text{tr} \mathbf A}{3 \varphi_0} \varepsilon, \quad
\varphi^{(1)} : = \varphi^{(0)} - R(\varphi^{(0)})/R'(\varphi^{(0)}), \quad
\varphi^{(2)} : = \varphi^{(1)} - R(\varphi^{(1)})/R'(\varphi^{(1)}).
\end{equation}
Finally, we put
\begin{equation}\label{2iterEBM}
{}^{n+1} \mathbf C_{\text{i}} =
\overline{ {}^{n+1} \mathbf C^{1/2}  \ \mathbf X(\varphi^{(2)}) \ {}^{n+1} \mathbf C^{1/2}}.
\end{equation}
We call this algorithm \emph{two iterations Euler backward method} (2IEBM).
The reason for choosing two iterations will be clear from the following sections.

\subsection{Properties of the IFEBM and 2IEBM}

Obviously, IFEBM and 2IEBM exactly preserve the geometric property \eqref{geopro}.
As shown in \cite{ShutovStability}, this allows one
to suppress the error accumulation when working with big time steps.
Moreover, the tensor $\mathbf{A}$ is positive definite. Thus, the
right-hand side of \eqref{ReducedEvolut5} is positive definite as well.
Taking \eqref{UpdateCi} and \eqref{UpdateCi_correct} into account, the solution
${}^{n+1} \mathbf C_{\text{i}}$ is also positive definite for both methods.

IFEBM is first order accurate (see Appendix C).
Note that IFEBM and 2IEBM exactly preserve the
w-invariance property (see Appendix D).
For the IFEBM, the solution
is well defined for all $\Delta t \geq 0$ and
${}^{n+1} \mathbf{C}_{\text{i}}$ it is a smooth function of the time step size $\Delta t$.
For $\Delta t \geq 0$, the solution ${}^{n+1} \mathbf{C}_{\text{i}}$ ranges
smoothly from ${}^{n} \mathbf{C}_{\text{i}}$ to ${}^{n+1} \overline{\mathbf C}$.
IFEBM is unconditionally stable since the solution remains
finite for arbitrary time steps.

\emph{Remark 2} \\
In the case of the neo-Hookean potential we have $c_{01} =0$.
Thus, $\varepsilon =0$, $\varphi = \varphi_0$, and $\mathbf{X} = \varphi^{-1} \mathbf{A}$.
Therefore, both methods are reduced to the the
explicit update formula from \cite{ShutovLandgraf}:
\begin{equation}\label{neoHookean}
c_{01} = 0 \quad  \Rightarrow  \quad
{}^{n+1} \mathbf C_{\text{i}} = \overline{{}^{n}\mathbf C_{\text{i}} +
\frac{\Delta t}{\eta}   c_{10} \overline{{}^{n+1} \mathbf C} }.
\end{equation}

\emph{Remark 3} \\
IFEBM is not
the only modification of EBM which ensures
the exact incompressibility.
Another modification of that kind was considered in \cite{Helm2}. Furthermore, in paper \cite{Vladimirov}
two other versions of the EBM were considered to enforce the inelastic
incompressibility by introducing the additional equation $\det \mathbf C_{\text{i}} = 1$. In contrast
to the IFEBM presented here, a local iterative
procedure was used in \cite{Helm2} and \cite{Vladimirov}.

\subsection{Explicit update formula in the Eulerian formulation}

The Eulerian formulation \eqref{SimMie82} of the flow rule is quite common.
Since its direct time discretization is not trivial (cf. \cite{Simo}), efficient time stepping
algorithms are needed.
In this subsection we transform the previously constructed IFEBM to build an iteration free
method on the current configuration. First, recall that for the
elastic left Cauchy-Green tensor ${\mathbf B}_{\text{e}} :=
{\mathbf F}_{\text{e}} {\mathbf F}^{\text{T}}_{\text{e}}$ the following relation holds
\begin{equation}\label{elsLCGT}
{\mathbf B}_{\text{e}} = \mathbf F \ \mathbf{C}^{-1}_{\text{i}} \ {\mathbf F}^{\text{T}}.
\end{equation}
For the generic time interval $(t_n, t_{n+1})$ we introduce
the so-called trial elastic left Cauchy-Green tensor ${}^{n+1} {\mathbf B}^{\text{trial}}_{\text{e}}$, by assuming that
$ \mathbf{C}_{\text{i}}$ remains constant during the step:
\begin{equation}\label{TrialBe}
{}^{n+1} {\mathbf B}^{\text{trial}}_{\text{e}} :=
{}^{n+1} {\mathbf F} \  {}^{n} \mathbf{C}^{-1}_{\text{i}} \ {}^{n+1} {\mathbf F}^{\text{T}}.
\end{equation}
Using the so-called relative deformation gradient
$\mathbf F_{\text{rel}} := {}^{n+1} \mathbf F \ {}^{n} \mathbf F^{-1}$ we have
\begin{equation}\label{TrialBerelat}
{}^{n+1} {\mathbf B}^{\text{trial}}_{\text{e}}  =
\mathbf F_{\text{rel}}  \ {}^{n} {\mathbf B}_{\text{e}} \ \mathbf F_{\text{rel}}^{\text{T}}.
\end{equation}
Next, we consider the polar decomposition
\begin{equation}\label{Polar}
{}^{n+1} \overline{\mathbf F} = {}^{n+1}\mathbf{R} \ ^{n+1} \overline{\mathbf C}^{1/2}, \quad
\text{where} \quad {}^{n+1}\mathbf{R} \ {}^{n+1}\mathbf{R}^{\text{T}} = \mathbf{1}, \quad \det({}^{n+1}\mathbf{R})=+1.
\end{equation}
Combining \eqref{elsLCGT} and $\eqref{Polar}_1$ with the definition of $\mathbf X$ given by $\eqref{ModEBM4}_1$, we have
\begin{equation}\label{Eulerian1}
\overline{{}^{n+1} {\mathbf B}^{-1}_{\text{e}}} = {}^{n+1} {\mathbf R} \ \mathbf{X} \  {}^{n+1} {\mathbf R}^{\text{T}}, \quad
(\overline{{}^{n+1} {\mathbf B}^{-1}_{\text{e}}})^2 = {}^{n+1} {\mathbf R} \ \mathbf{X}^2 \  {}^{n+1} {\mathbf R}^{\text{T}}.
\end{equation}
Pre-multiplying \eqref{ModEBM3} with ${}^{n+1} {\mathbf R}$ and post-multiplying it with ${}^{n+1} {\mathbf R}^{\text{T}}$,
we obtain the following quadratic equation with respect to unknown $\overline{{}^{n+1} {\mathbf B}^{-1}_{\text{e}}}$
\begin{equation}\label{Eulerian2}
\varphi \ \overline{{}^{n+1} {\mathbf B}^{-1}_{\text{e}}} = {}^{n+1} {\mathbf R} \ \mathbf{A} \ {}^{n+1} {\mathbf R}^{\text{T}}
 - \varepsilon \ (\overline{{}^{n+1} {\mathbf B}^{-1}_{\text{e}}})^{2}.
\end{equation}
To simplify this relation we introduce $\tilde{\mathbf{A}} := {}^{n+1} {\mathbf R} \
\mathbf{A} \ {}^{n+1}{\mathbf R}^{\text{T}}$; then
\begin{equation}\label{Eulerian3}
\tilde{\mathbf{A}} = {}^{n+1} {\mathbf R} \ \overline{{}^{n+1} \mathbf C^{-1/2}} \big(
 {}^{n}\mathbf C_{\text{i}} + \frac{\Delta t}{\eta} c_{10}
 \overline{{}^{n+1} \mathbf C} ) \overline{{}^{n+1}
 \mathbf C^{-1/2}} \ {}^{n+1} {\mathbf R}^{\text{T}} =
 \overline{({}^{n+1} {\mathbf B}^{\text{trial}}_{\text{e}} )^{-1}} +
 \frac{\Delta t}{\eta} c_{10} \mathbf{1}.
\end{equation}
Thus, we arrive at the following problem
\begin{equation}\label{Eulerian4}
\varphi \ \overline{{}^{n+1} {\mathbf B}^{-1}_{\text{e}}} =  \tilde{\mathbf{A}}
 - \varepsilon \ (\overline{{}^{n+1} {\mathbf B}^{-1}_{\text{e}}})^{2}, \quad
 \det (\overline{{}^{n+1} {\mathbf B}^{-1}_{\text{e}}}) =1.
\end{equation}
In analogy to \eqref{ReducedEvolut5notCont}, the following closed form solution is valid
\begin{equation}\label{Eulerian5}
\overline{{}^{n+1} {\mathbf B}^{-1}_{\text{e}}} = 2  \tilde{\mathbf{A}} \Big[  \Big( \varphi^2 \textbf{1} +
4 \varepsilon  \tilde{\mathbf{A}} \Big)^{1/2} + \varphi \textbf{1} \Big]^{-1}.
\end{equation}
Taking into account that $\det \tilde{\mathbf{A}} = \det \mathbf{A}$ and $\text{tr} \tilde{\mathbf{A}} = \text{tr} \mathbf{A}$, we have
the spatial analog of \eqref{EstimOfz}:
\begin{equation}\label{Eulerian6}
\varphi = \varphi_0 - \frac{\text{tr} \tilde{\mathbf{A}}}{3 \varphi_0} \varepsilon + O ( \varepsilon^2),
\quad \text{where} \
\varphi_0 :=  \big(\det \tilde{\mathbf{A}} \big)^{1/3}.
\end{equation}
Since \eqref{Eulerian6} is not exact, an additional correction is needed:
\begin{equation}\label{Eulerian5corrected}
\overline{{}^{n+1} {\mathbf B}^{-1}_{\text{e}}} = \overline{2  \tilde{\mathbf{A}} \Big[  \Big( \varphi^2 \textbf{1} +
4 \varepsilon  \tilde{\mathbf{A}} \Big)^{1/2} + \varphi \textbf{1} \Big]^{-1}}.
\end{equation}

The \emph{explicit} update procedure for the evolution equation \eqref{SimMie8} is summarized in table \ref{tab2}.
Within a time step, the explicit procedure \eqref{Eulerian6}-\eqref{Eulerian5corrected} predicts
the same stress response as the previously reported Lagrangian procedure.
Note that due to the elastic isotropy
${}^{n+1} {\mathbf B}_{\text{e}}$ and ${}^{n+1} {\mathbf B}^{\text{trial}}_{\text{e}}$ are coaxial (cf. \cite{Simo}).

\begin{table}[h!]
\caption{IFEBM on the current configuration}
\begin{center}
\begin{tabular}{|l |}
\hline
input: ${}^{n+1} \mathbf F$, ${}^{n} \mathbf F$, $\overline{{}^{n} \mathbf B_{\text{e}}^{-1}}$  \\
output: $\overline{{}^{n+1} \mathbf B_{\text{e}}^{-1}}$, ${}^{n+1} \mathbf S$  \\
1: $\overline{\mathbf F_{\text{rel}}} = \overline{{}^{n+1} \mathbf F \ {}^{n} \mathbf F^{-1}}$ \\
2: $\overline{({}^{n+1} {\mathbf B}^{\text{trial}}_{\text{e}})^{-1}}  =
(\overline{\mathbf F_{\text{rel}}})^{-\text{T}}  \ \overline{{}^{n}
{\mathbf B}^{-1}_{\text{e}}} \ (\overline{\mathbf F_{\text{rel}}})^{-1}$    \\
3: $\tilde{\mathbf{A}} = \overline{({}^{n+1} {\mathbf B}^{\text{trial}}_{\text{e}} )^{-1}} + \frac{\Delta t}{\eta} c_{10} \mathbf{1}$  \\
4: $\varepsilon = c_{01} \frac{\Delta t}{\eta}$    \\
5: $\varphi_0 =  \big(\det \tilde{\mathbf{A}} \big)^{1/3}$ \\
6: $\varphi = \varphi_0 - \frac{\text{tr} \tilde{\mathbf{A}}}{3 \varphi_0} \varepsilon$  \\
7: $\overline{{}^{n+1} {\mathbf B}^{-1}_{\text{e}}} = \overline{2  \tilde{\mathbf{A}} \Big[  \Big( \varphi^2 \textbf{1} +
4 \varepsilon  \tilde{\mathbf{A}} \Big)^{1/2} + \varphi \textbf{1} \Big]^{-1}}$  \\
8: $^{n+1} \mathbf S = c_{10} (\overline{^{n+1} \mathbf{B}_{\text{e}}})^{\text{D}}
-c_{01} (\overline{^{n+1} \mathbf{B}^{-1}_{\text{e}}})^{\text{D}} $    \\  \hline
\end{tabular} \\
\end{center}
\label{tab2}
\end{table}

\section{Numerical tests}

\subsection{Non-proportional loading}

All quantities in this subsection are non-dimensional.
Assume the following local deformation history
\begin{equation}\label{loaprog0}
\mathbf F (t) = \overline{\mathbf F^{\prime} (t)} \quad \text{for} \quad t \in [1, 3],
\end{equation}
where $\mathbf F^{\prime} (t)$ is a piecewise linear function
\begin{equation}\label{loaprog}
\mathbf F^{\prime} (t) :=
\begin{cases}
    (1 - t) \mathbf F_1  + t \ \mathbf F_2 \quad \quad  \ \  \text{if} \ t \in [0,1] \\
    (2 - t) \mathbf F_2  + (t-1) \mathbf F_3 \quad \text{if} \ t \in (1,2] \\
     (3 - t) \mathbf F_3  + (t-2) \mathbf F_4 \quad \text{if} \ t \in (2,3]
\end{cases},
\end{equation}
\begin{equation}\label{loaprog2}
\mathbf F_1 :=\mathbf 1, \quad
\mathbf F_2 := 2 \mathbf{e}_{1} \otimes \mathbf{e}_{1} +
\frac{1}{\sqrt2} (\mathbf{e}_{2} \otimes \mathbf{e}_{2} + \mathbf{e}_{3} \otimes \mathbf{e}_{3}),
\end{equation}
\begin{equation*}\label{loaprog222}
\mathbf F_3 := \mathbf 1 + \mathbf{e}_{1} \otimes \mathbf{e}_{2}, \quad
\mathbf F_4 := 2 \mathbf{e}_{2} \otimes \mathbf{e}_{2} +
\frac{1}{\sqrt2} (\mathbf{e}_{1} \otimes \mathbf{e}_{1} + \mathbf{e}_{3} \otimes \mathbf{e}_{3}).
\end{equation*}
This  non-proportional loading programm is strain-driven and
volume-preserving; abrupt changes of the loading direction appear
at $t=1$ and $t=2$. The initial condition is
$\mathbf C_{\text{i}}|_{t=0} =  \mathbf 1$; it means
that the reference configuration is stress free at $t=0$.

The numerical solution of the initial value problem \eqref{puba44}, \eqref{initCond}
obtained with very small time steps
will be referred to as the exact solution. Let ${\mathbf S}^{\text{exact}}(t)$ be the corresponding history of the
Kirchhoff stress.
Next, large time steps $\Delta t$ are utilized to test the accuracy of the numerical algorithms,
 thus leading to large strain increments. The corresponding computed
Kirchhoff stress is denoted by ${\mathbf S}^{\text{num}}(t)$, and the corresponding
error in stress prediction is given by the Frobenius norm
\begin{equation}\label{StressError}
\text{Error}(t) : = \| {\mathbf S}^{\text{exact}}(t) - {\mathbf S}^{\text{num}}(t) \|.
\end{equation}
The following constants of the Mooney-Rivling elasticity are used: $c_{10} = c_{01} = 1$.
The proposed IFEBM is compared with a modified Euler backward method
(MEBM) and the classical exponential method (EM). A short summary of these methods is provided in Appendix E.
The error is plotted versus time in Figure \ref{fig1} for $\Delta t = 0.1$ and in Figure \ref{fig2} for $\Delta t = 0.05$.
The figures show that IFEBM, MEBM, and EM are equivalent in terms of accuracy.
The difference between IFEBM and MEBM is much smaller than the difference between MEBM and EM.\footnote{
The two iteration Euler backward method (2IEBM) is tested as well. The corresponding error curve
practically merges with the error curve for MEBM. Therefore, it is not depicted in Figures \ref{fig1} and \ref{fig2}.}
Since the considered methods are first order accurate, the error
for $\Delta t=0.05$ is approximately two times smaller than for $\Delta t=0.1$.

The Newton-Raphson procedure for finding the numerical
solution, pertaining to MEBM and EM, may diverge, if the solution from the previous time step
is taken as the initial approximation. A substepping can be used to resolve this problem.
On the other hand, the solution for the novel
IFEBM is given in a closed form. Thus, the new method is a
priory free from any convergence problems, which makes it more robust.

\begin{figure}[h!]\centering
\scalebox{0.85}{\includegraphics{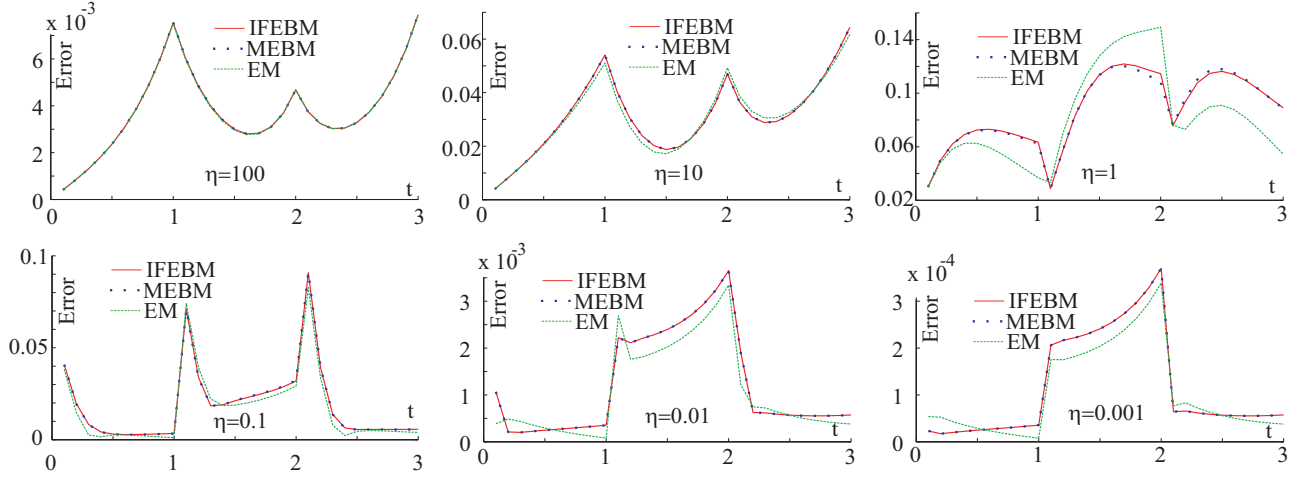}}
\caption{Error graphs for different integrations methods; $\Delta t = 0.1$. \label{fig1}}
\end{figure}

\begin{figure}[h!]\centering
\scalebox{0.85}{\includegraphics{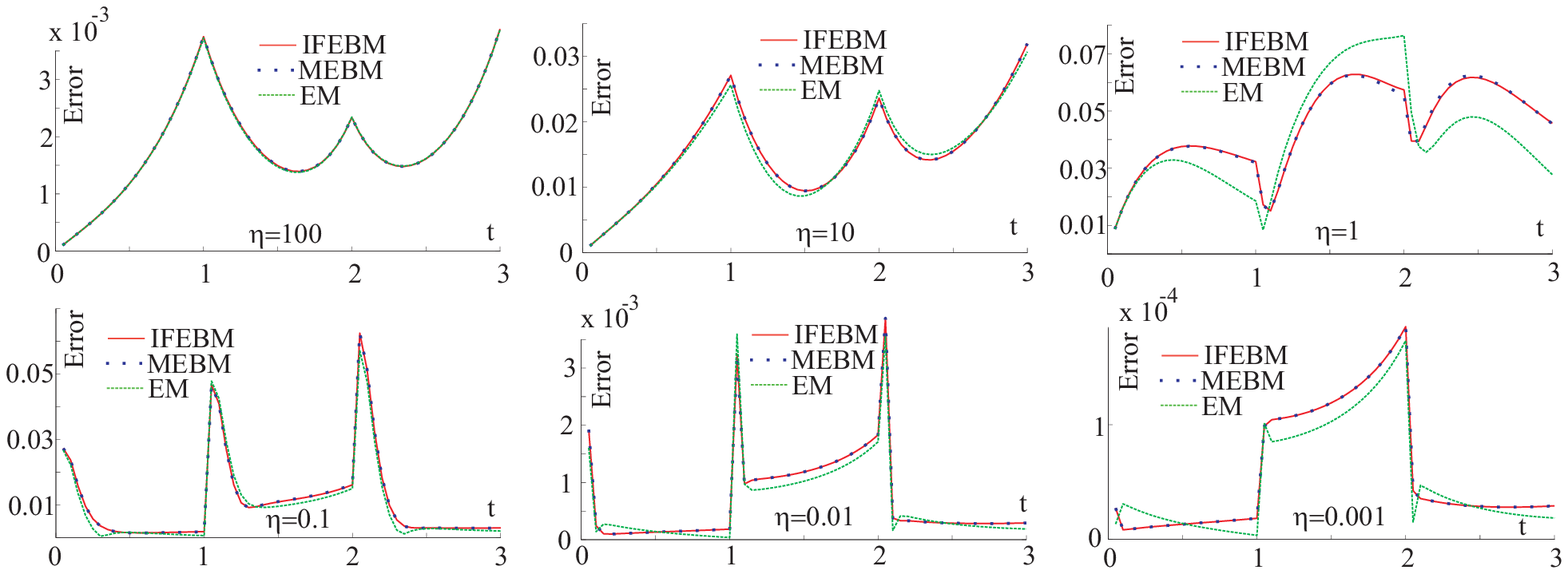}}
\caption{Error graphs for different integrations methods;  $\Delta t = 0.05$. \label{fig2}}
\end{figure}

Another  aspect is the symmetry of the consistent tangent operator
$ \frac{\displaystyle \partial \ {}^{n+1} \tilde{\mathbf{T}}}{\displaystyle \partial \ {}^{n+1} \mathbf{C}}$.
In the current study the tangent is computed by numerical differentiation
using the central difference scheme, which provides
at least ten-digit accuracy.
For computations, the relevant tensors are represented by the
6-vectors as follows (cf. the Voigt notation)
\begin{equation}\label{VectorTensorT}
\overrightarrow{\mathbf{T}}:=({}^{n+1} \tilde{\mathbf{T}}_{11},
{}^{n+1} \tilde{\mathbf{T}}_{22}, {}^{n+1} \tilde{\mathbf{T}}_{33}, {}^{n+1} \tilde{\mathbf{T}}_{12},
{}^{n+1} \tilde{\mathbf{T}}_{13}, {}^{n+1} \tilde{\mathbf{T}}_{23}),
\end{equation}
\begin{equation}\label{VectorTensorC}
 \overrightarrow{\mathbf{C}}:=({}^{n+1} \mathbf{C}_{11}, {}^{n+1} \mathbf{C}_{22}, {}^{n+1} \mathbf{C}_{33},
2 \ {}^{n+1} \mathbf{C}_{12}, 2 \ {}^{n+1} \mathbf{C}_{13}, 2 \ {}^{n+1} \mathbf{C}_{23}).
\end{equation}
The deviation of the tangent from the symmetry is measured by the following normalized quantity
\begin{equation}\label{SymDeviation}
\text{Deviation} := \max_{t \in [0,3]} \Big\| \frac{\displaystyle \partial \ \overrightarrow{\mathbf{T}}}{\displaystyle
\partial \  \overrightarrow{\mathbf{C}}} - \Big( \frac{\displaystyle \partial \ \overrightarrow{\mathbf{T}}}{\displaystyle
\partial \  \overrightarrow{\mathbf{C}}}  \Big)^{\text{T}} \Big\| \Big/ \max_{t \in [0,3]}  \Big\|  \frac{\displaystyle \partial \
\overrightarrow{\mathbf{T}}}{\displaystyle \partial \  \overrightarrow{\mathbf{C}}} \Big\|, \quad
\| M \| : = \Big( \sum^6_{i,j=1} M^2_{ij} \Big)^{1/2}.
\end{equation}
The computed deviation of the consistent tangent from the symmetry
is summarized in Table \ref{tab_IFEBM} for IFEBM
and in Table \ref{tab_2IEBM} for 2IEBM. As can be seen from the table,
IFEBM provides a tangent which is close to a symmetric one. For the
2IEBM the deviation from the symmetry is barely detectable.
Thus, 2IEBM can be used in applications, where the symmetry of the tangent becomes crucial.

\begin{table}[h]
\caption{IFEBM: deviation of the consistent tangent operator from the symmetry.}
\begin{center}
\begin{tabular}{| l|  l |  l  | l | l | l | l |}
\hline
            & $\eta = 100$       &  $\eta = 10$   &   $\eta = 1$ & $\eta = 0.1$ & $\eta = 0.01$ & $\eta = 0.001$     \\ \hline
 $\Delta t=0.1$  &  $< 10^{-9}$       &   $8.5 \cdot 10^{-7}$      &   $1.8 \cdot 10^{-4}$  &  $8.5 \cdot 10^{-5}$   &  $6.0 \cdot 10^{-7}$   & $< 10^{-9}$     \\ \hline
 $\Delta t=0.05$ &   $< 10^{-9}$      &   $1.1 \cdot 10^{-7}$     &  $3.0 \cdot 10^{-5}$   &  $2.7 \cdot 10^{-5}$   &  $4.9 \cdot 10^{-7}$  & $< 10^{-9}$     \\ \hline
\end{tabular}
\end{center}
\label{tab_IFEBM}
\end{table}

\begin{table}[h]
\caption{2IEBM: deviation of the consistent tangent operator from the symmetry.}
\begin{center}
\begin{tabular}{| l|  l |  l  | l | l | l | l |}
\hline
            & $\eta = 100$       &  $\eta = 10$   &   $\eta = 1$ & $\eta = 0.1$ & $\eta = 0.01$ & $\eta = 0.001$     \\ \hline
 $\Delta t=0.1$  &  $< 10^{-9}$       &   $< 10^{-9}$      &   $< 10^{-9}$  &  $1.2 \cdot 10^{-9}$   &  $< 10^{-9}$   & $< 10^{-9}$     \\ \hline
 $\Delta t=0.05$ &   $< 10^{-9}$      &   $< 10^{-9}$     &  $< 10^{-9}$   &  $< 10^{-9}$   &  $< 10^{-9}$  & $< 10^{-9}$     \\ \hline
\end{tabular}
\end{center}
\label{tab_2IEBM}
\end{table}

\emph{Remark 4} \\
A popular way to enforce the symmetry of the tangent operator
exploits the variational nature of the underlying constitutive equations
\cite{Fancello, BleierMosler, Vassoler}. However, iterative procedures are used in these references.

\emph{Remark 5} \\
In a number of FEM applications, the symmetry of the tangent operator is not important.
These applications include globally explicit FEM procedures (the stiffness
matrix is not assembled at all),
geometrically nonlinear computations with follower loads or
computations with nonlinear kinematic hardening
(the stiffness matrix is a priory non symmetric).
Moreover, numerical tests from Section 4.3 show that the global
convergence of the FEM procedure with artificially symmetrized
tangent is still very good.

\subsection{Uniaxial compression tests}

In this subsection we consider a composite model
with one equilibrium branch and four Maxwell branches connected in parallel.
The free energy of the equilibrium branch is given by
\begin{equation}\label{Equilibrium}
\rho_{\scriptscriptstyle \text{R}} \psi_{\text{eq}}(\mathbf{C}) =
\frac{c^{(\text{eq})}_{10}}{2} \big( \text{tr} \overline{\mathbf{C}} - 3 \big) +
\frac{c^{(\text{eq})}_{01}}{2} \big( \text{tr} \overline{\mathbf{C}} - 3 \big)
+ \frac{k}{50}\big( (\text{det} \mathbf{C})^{5/2} +
(\text{det} \mathbf{C})^{-5/2} -2 \big).
\end{equation}
The volumetric part appearing in \eqref{Equilibrium} is taken from \cite{HartmannNeff};
$k$ is the bulk modulus.
For the m-th Maxwell branch we put
\begin{equation}\label{specLike} 
\rho_{\scriptscriptstyle \text{R}}  \psi_{\text{ov},m}(\hat{\mathbf{C}}^{(m)}_{\text{e}})=
\frac{c^{(m)}_{10}}{2} \big( \text{tr} \overline{\hat{\mathbf{C}}^{(m)}_{\text{e}}} - 3 \big) +
\frac{c^{(m)}_{01}}{2} \big( \text{tr} \overline{(\hat{\mathbf{C}}^{(m)}_{\text{e}})^{-1}} - 3 \big), \ m=1,2,3,4.
\end{equation}
For the composite model, the 2nd Piola-Kirchhoff stress equals
\begin{equation}\label{SumStresses}
\tilde{\mathbf{T}}
= \tilde{\mathbf{T}}_{\text{eq}} + \sum\limits_{m=1}^{4} \tilde{\mathbf{T}}_{\text{ov},m}, \quad \text{where}
\end{equation}
\begin{equation}\label{EquilibriumStress}
\tilde{\mathbf{T}}_{\text{eq}}
= \mathbf C^{-1} (c^{(\text{eq})}_{10} \overline{\mathbf C}
- c^{(\text{eq})}_{01} \overline{\mathbf C^{-1}} )^{\text{D}} + \frac{\displaystyle k}{\displaystyle 10} \
\big( (\text{det} \mathbf C)^{5/2}-(\text{det} \mathbf C)^{-5/2} \big) \
\mathbf C^{-1},
\end{equation}
\begin{equation}\label{2PKdover}
\tilde{\mathbf T}_{\text{ov},m}  = \mathbf C^{-1} \big(c^{(m)}_{10} \overline{\mathbf C} (\mathbf C^{(m)}_{\text{i}})^{-1}
- c^{(m)}_{01} \mathbf C^{(m)}_{\text{i}} \overline{\mathbf C^{-1}} \big)^{\text{D}}.
\end{equation}
Analogously to \eqref{puba44}, the behaviour of each Maxwell branch is described by
\begin{equation}\label{EvolMaxwe}
\dot{\mathbf C}^{(m)}_{\text{i}} = \frac{1}{\eta^{(m)}}  \big( c^{(m)}_{10} \overline{\mathbf C}
(\mathbf C_{\text{i}}^{(m)})^{-1} - c^{(m)}_{01}
\mathbf C^{(m)}_{\text{i}}  \overline{\mathbf C^{-1}} \big)^{\text{D}} \mathbf C^{(m)}_{\text{i}}, \quad
\mathbf C^{(m)}_{\text{i}}|_{t=0} = \mathbf{1}.
\end{equation}
The material parameters corresponding to a cartilaginous temporomandibular joint are
taken from \cite{Koolstra}; they are summarized in Table \ref{BioMat}.

\begin{table}[h]
\caption{Material parameters of a cartilaginous temporomandibular joint \cite{Koolstra}.}
\begin{center}
\begin{tabular}{| l| l | l | l | l | l |}
\hline
            & m =1            &  m = 2    &   m =3 & m=4 & equilibrium      \\ \hline
 $c_{10}$ [MPa]  &   0.25     &   0.25    &  0.36    &  1.25   &  0.2     \\ \hline
 $c_{01}$ [MPa]  &   0.25     &   0.25    &  0.36    &  1.25   &  0.2     \\ \hline
 $\eta$ [MPa s]  &   25.0     &   5.0     &  0.144   &  0.005  &  $\infty$     \\ \hline
\end{tabular}
\end{center}
\label{BioMat}
\end{table}

\begin{figure}[t]\centering
\scalebox{0.85}{\includegraphics{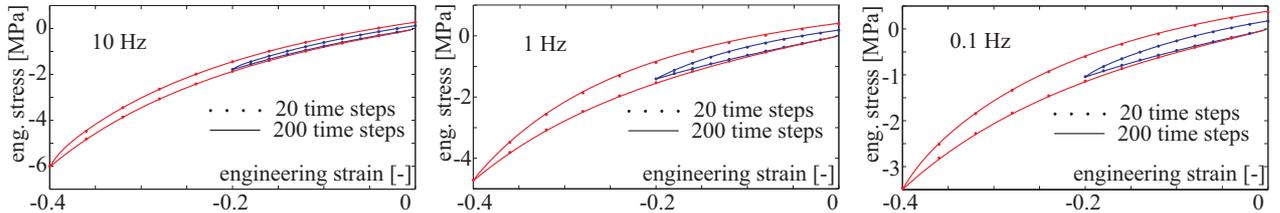}}
\caption{IFEBM-based simulations of a non-monotonic
compression test with different strain amplitudes and frequencies. \label{UniaxTenssCompress}}
\end{figure}

In this subsection we visualize the stress response of the joint tissue
 subjected to a non-monotonic volume-preserving uniaxial loading.
Here, the material is assumed to be incompressible.\footnote{Formally, this corresponds to
$k \rightarrow  \infty$ in ansatz \eqref{Equilibrium}.}
The strain-controlled loading is given by
\begin{equation}\label{IncompreTensCompr}
\mathbf F(t) = (1+\varepsilon(t)) \mathbf{e}_{1} \otimes \mathbf{e}_{1} + (1+\varepsilon(t))^{-1/2} (
\mathbf{e}_{2} \otimes \mathbf{e}_{2} + \mathbf{e}_{3} \otimes \mathbf{e}_{3}),
\end{equation}
where $\varepsilon$ stands for the prescribed engineering strain.
Three different oscillation frequencies are analyzed (10 Hz, 1 Hz, and 0.1 Hz).
For each frequency, two different tests are simulated: one test with the
strain amplitude 0.2 and another one with the strain amplitude 0.4.
Within each test, the absolute strain rate is constant: $|\dot{\varepsilon}| = const$.
The simulation results are shown in Figure \ref{UniaxTenssCompress}.
The figure reveals that for different loading frequencies IFEBM is sufficiently
accurate, even when working with moderate time steps.

\subsection{Solution of a boundary value problem}

In order to demonstrate the applicability of the IFEBM,
a representative initial boundary value problem is solved here.
The composite material model
\eqref{EquilibriumStress}--\eqref{EvolMaxwe}
 introduced in the previous subsection is implemented into the commercial FEM code MSC.MARC via
the Hypela2 user interface. The material parameters correspond to the
cartilaginous temporomandibular joint (cf. Table \ref{BioMat})
with the bulk modulus $k = 20$ MPa.

\begin{figure}[h!]\centering
\scalebox{0.85}{\includegraphics{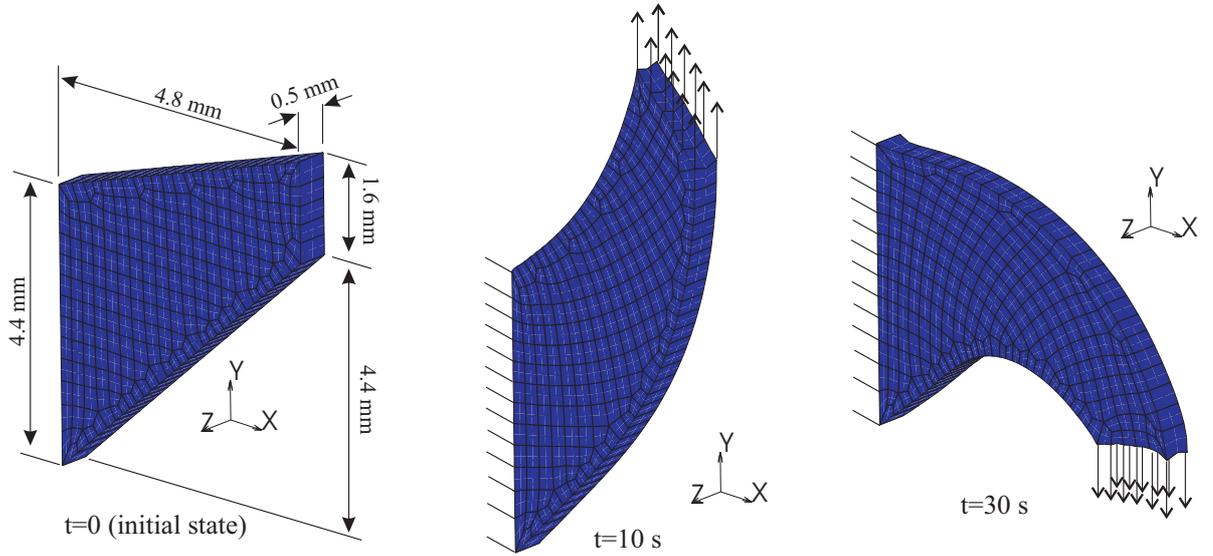}}
\caption{Force-controlled loading of the Cook membrane. Left: specimen's geometry and dimensions. Middle and right:
boundary conditions and deformed configurations at different instances of time. \label{FEMCooksMembraine}}
\end{figure}

A non-monotonic force-controlled loading is applied to the Cook membrane
(Figure \ref{FEMCooksMembraine}).\footnote{The preference is given to the force-controlled loading
here since such problems are more difficult to solve than the strain-controlled ones.}
The membrane is discretized using 213 elements with quadratic approximation of
geometry and displacements; full integration is used.
Let $T_{\text{total}}$ be the overall process time.
For $t \in [0, T_{\text{total}}/3]$ the applied load
increases from zero to the maximum of $0.27$ N;
after that the load linearly decreases to $-0.27$ N.
Three different loading rates are analyzed, leading to
different process durations: $T_{\text{total}} = 0.3$ s,
$T_{\text{total}} = 3$ s, and $T_{\text{total}} = 30$ s.
The deformed shapes are shown in Figure \ref{FEMCooksMembraine}
for $T_{\text{total}} = 30$ at different time instances.

Constant time steps are implemented in each simulation; simulations with small and
large time steps are carried out.
A strict convergence criterion is used at each time step: for convergence, the relative
error in force and displacements should not exceed 0.01.\footnote{The MSC.MARC settings are as follows:
relative force tolerance = 0.01, relative displacement tolerance = 0.01.}
Two types of computations are carried out: with a general non-symmetric stiffness matrix and with an
artificially symmetrized stiffness matrix. The symmetrized version allows us to use a
symmetric matrix solver, which is more efficient.
Although the consistent tangent operator for IFEBM is not
exactly symmetric, the artificial symmetrization did not affect the number
of the global equilibrium iterations: MSC.MARC required the same number
of iterations,
both with symmetric and non-symmetric matrix solvers.
The required number of iterations is listed in Table \ref{NumOfIter}.
More iterations are needed for the slow process.

\begin{table}[h]
\caption{Number of global equilibrium iterations.}
\begin{center}
\begin{tabular}{| l| l | l | l |}
\hline
            & $T_{\text{total}} = 0.3$ s &  $T_{\text{total}} = 3$ s &  $T_{\text{total}} = 30$ s      \\ \hline
 150 steps  &   449     &   448    &  480    \\ \hline
 15 steps   &   58      &   70     &  131    \\ \hline
\end{tabular}
\end{center}
\label{NumOfIter}
\end{table}

\begin{figure}[h!]\centering
\scalebox{0.85}{\includegraphics{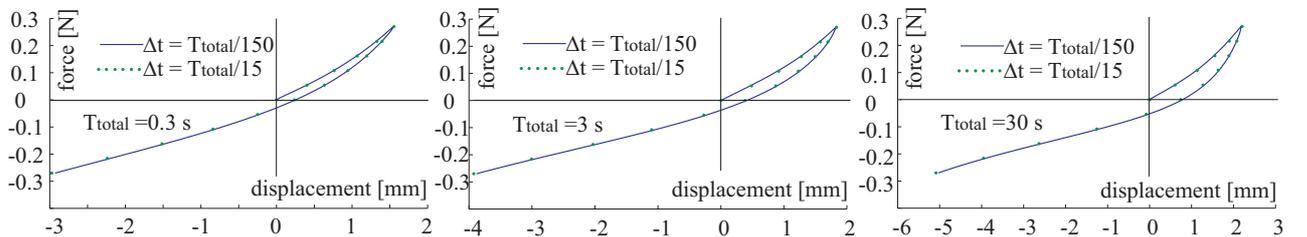}}
\caption{Simulated force-displacement curves for three processes
with different durations. \label{ForceDisplacmentCurves}}
\end{figure}

The simulated force-displacement curves are shown in Figure \ref{ForceDisplacmentCurves}.
Although the stress response is history-dependent, simulations with only 15 steps provide
sufficiently accurate results. The maximum displacement increases with decreasing loading rate.
Thus, the large number of the global equilibrium iterations for the slow process can be explained
by pronounced geometric nonlinearities.

\section{Discussion and conclusion}

A popular version of the finite strain Maxwell fluid is considered, which is based
on the multiplicative decomposition of the deformation gradient, combined
with hyperelastic relations between stresses and elastic strains.
A new iteration free method is proposed for the implicit time stepping,
provided that the elastic potential is of the Mooney-Rivlin type.

The following properties of the exact solution are inherited by the proposed IFEBM:
\begin{itemize}
\item[i] ${}^{n+1} \mathbf{C}_{\text{i}}$ is symmetric and unimodular: ${}^{n+1} \mathbf{C}_{\text{i}}
\in Sym$, $\det({}^{n+1} \mathbf{C}_{\text{i}}) = 1$,
\item[ii] ${}^{n+1} \mathbf{C}_{\text{i}}$ is positive definite: ${}^{n+1}\mathbf{C}_{\text{i}} > 0$,
\item[iii] ${}^{n+1}\mathbf{C}_{\text{i}}$ is a
smooth function of $\Delta t$,
\item[iv] complete stress relaxation: ${}^{n+1} \mathbf{C}_{\text{i}} \rightarrow {}^{n+1}
\overline{\mathbf{C}}$ as $\Delta t \rightarrow \infty$,
\item[v] w-invariance under
volume-preserving changes of the reference configuration.
\end{itemize}
For IFEBM, the consistent tangent operator is nearly symmetric.
A slight modification, called 2IEBM, allows one to obtain a tangent
operator which is even closer
to the set of symmetric fourth-rank tensors.
In terms of accuracy, the IFEBM is equivalent to the existing methods like the
Euler backward method with exact incompressibility of the flow (MEBM) and
the exponential method (EM), but the new method is superior in
terms of robustness and computational
efficiency.

Since the considered
multiplicative Maxwell fluid is widespread in the phenomenological material
modelling, the new method can become a method
of choice in many applications, especially in those that require increased
robustness and efficiency.

\textit{Acknowledgements}
This research was supported by the Russian Science Foundation (project number 15-11-20013).

\section*{Appendix A: alternative Eulerian formulation of the model}

Let us consider an alternative formulation of the finite strain Maxwell model from Section 2.
For simplicity, only incompressible behaviour is analyzed here.
Consider the velocity gradient $\mathbf{L}$, the continuum spin $\mathbf{W}$ and the strain rate $\mathbf{D}$
\begin{equation}\label{Kinemat}
\mathbf{L} := \dot{\mathbf{F}} \mathbf{F}^{-1}, \quad \mathbf{W}:=(\mathbf{L}  - \mathbf{L}^{\text{T}})/2, \quad
\mathbf{D}:=(\mathbf{L}  + \mathbf{L}^{\text{T}})/2.
\end{equation}
The formulation
is based on the additive decomposition of the strain rate tensor $\mathbf{D}$ into the elastic part $\mathbf{D}_{\text{e}}$
and the inelastic part $\mathbf{D}_{\text{i}}$
\begin{equation}\label{Additive}
\mathbf{D} = \mathbf{D}_{\text{e}} + \mathbf{D}_{\text{i}}.
\end{equation}
Following \cite{Donner}, assume that
$\mathbf{X}_\text{e}$ is a symmetric positive definite tensor with $\det(\mathbf{X}_\text{e}) = 1$;
$\mathbf{X}_\text{e}$ operates on the current configuration.
Let the free energy per unit mass be given by the isotropic function
 $\psi(\mathbf{X}_\text{e})$. Assume that the
Kirchhoff stress $\mathbf{S}$ is computed through
\begin{equation}\label{Murnag}
\mathbf{S} =  2 \rho_{\scriptscriptstyle \text{R}}
\frac{\displaystyle \partial \psi(\mathbf{X}_\text{e})}
{\displaystyle \partial \mathbf{X}_\text{e}} \mathbf{X}_\text{e}.
\end{equation}
The flow rule is postulated on the current configuration
\begin{equation}\label{FlowCurrent}
\mathbf{D}_{\text{i}} = \frac{1}{2 \eta} \mathbf{S}^{\text{D}}.
\end{equation}
Moreover, the evolution of $\mathbf{X}_\text{e}$ is described by (cf. \cite{Donner})
\begin{equation}\label{EvolXe}
\mathbf{X}^{ZJ}_\text{e} = \mathbf{D}_{\text{e}} \ \mathbf{X}_{\text{e}} + \mathbf{X}_{\text{e}} \ \mathbf{D}_{\text{e}},
\quad \text{where} \quad \mathbf{X}^{ZJ}_\text{e}: = \dot{\mathbf{X}}_{\text{e}}+ \mathbf{X}_{\text{e}} \ \mathbf{W} - \mathbf{W} \ \mathbf{X}_{\text{e}}.
\end{equation}
Let us show that this model is equivalent to \eqref{SimMie82}-\eqref{Murnag0}.
First, combining \eqref{Additive}  and \eqref{EvolXe}, we arrive at
\begin{equation}\label{AppA1}
\dot{\mathbf{X}}_{\text{e}} + \mathbf{X}_{\text{e}} \ \mathbf{W} - \mathbf{W} \ \mathbf{X}_{\text{e}} =
(\mathbf{D} - \mathbf{D}_{\text{i}}) \mathbf{X}_{\text{e}} + \mathbf{X}_{\text{e}} (\mathbf{D} - \mathbf{D}_{\text{i}}).
\end{equation}
Next, substituting \eqref{FlowCurrent} into this, we obtain after some rearrangements
\begin{equation}\label{AppA2}
\dot{\mathbf{X}}_{\text{e}} - \mathbf{X}_{\text{e}} \ \mathbf{L}^{\text{T}} - \mathbf{L} \ \mathbf{X}_{\text{e}} =
 - \frac{1}{2 \eta} (\mathbf{S}^{\text{D}} \mathbf{X}_{\text{e}} + \mathbf{X}_{\text{e}} \mathbf{S}^{\text{D}}).
\end{equation}
Relation \eqref{Murnag} implies that $\mathbf{X}_{\text{e}}$ and $\mathbf{S}^{\text{D}}$ commute. Post-multiplying
both sides of \eqref{AppA2} with $\mathbf{X}^{-1}_{\text{e}}$ and
recalling the definition of the Lie derivative \eqref{OldrCovar}, we obtain
\begin{equation}\label{AppA3}
\text{\calligra{L}}_{v}  (\mathbf{X}_{\text{e}}) \mathbf{X}^{-1}_{\text{e}} =  - \frac{1}{\eta} \mathbf{S}^{\text{D}}.
\end{equation}
It remains to note that \eqref{Additive} and \eqref{AppA3} become identical to \eqref{SimMie82}-\eqref{Murnag0},
if $\mathbf{X}_{\text{e}}$ is formally replaced by $\mathbf{B}_{\text{e}}$.
Thus, we are dealing with two different Eulerian formulations of the same model.

\section*{Appendix B: estimation of $\varphi$}

Here we estimate the unknown $\varphi$; the derivation is similar to that presented in Appendix C of \cite{Shutov2016}.
The equation $\varphi \ \mathbf X = \mathbf A -  \varepsilon \mathbf X^2$ yields $\mathbf X$
as an implicit function of $\varphi$ and $\varepsilon$.
Its expansion in Taylor series for small  $\varepsilon$ is as follows
\begin{equation}\label{AppendixB3}
\mathbf X = \tilde{\mathbf X}(\varphi,\varepsilon) =
\frac{1}{\varphi} \mathbf A - \frac{\varepsilon}{\varphi^3} {\mathbf A}^2 + O(\varepsilon^2), \quad
\mathbf X|_{\varepsilon =0} = \frac{1}{\varphi} \mathbf A.
\end{equation}
The unknown $\varphi$ is estimated using the incompressibility relation $\det (\textbf{X}) = 1$, which yields
\begin{equation}\label{AppendixB4}
\varphi = \tilde{\varphi} (\varepsilon), \quad \varphi_0 := \tilde{\varphi} (0)
= (\det \mathbf A )^{1/3}.
\end{equation}
Using the implicit function theorem, we have
\begin{equation}\label{AppendixB5}
\frac{d \tilde{\varphi} (\varepsilon)}{d \varepsilon}|_{\varepsilon =0} =
- \frac{\partial \det \tilde{\mathbf X}(\varphi,\varepsilon)}{\partial
\varepsilon}|_{\varphi=\varphi_0, \varepsilon =0}
\ \Big(\frac{\partial \det \tilde{\mathbf X}(\varphi,\varepsilon)}
{\partial \varphi}|_{\varphi=\varphi_0, \varepsilon =0}\Big)^{-1}.
\end{equation}
Differentiating expansion $\eqref{AppendixB3}_1$, we obtain
\begin{equation}\label{AppendixB52}
\frac{\partial \det \tilde{\mathbf X}(\varphi,\varepsilon)}{\partial \varepsilon}|_{\varphi=\varphi_0, \varepsilon =0} =
\det \tilde{\mathbf X}(\varphi_0,0) \ (\tilde{\mathbf X}(\varphi_0,0))^{-1} : \Big(-\frac{\mathbf A^2}{\varphi_0^3}\Big),
\end{equation}
\begin{equation}\label{AppendixB53}
\frac{\partial \det \tilde{\mathbf X}(\varphi,\varepsilon)}{\partial \varphi}|_{\varphi=\varphi_0, \varepsilon =0} =
\det \tilde{\mathbf X}(\varphi_0,0) \ (\tilde{\mathbf X}(\varphi_0,0))^{-1} : \Big(- \frac{\mathbf A}{\varphi_0^2}\Big).
\end{equation}
Substituting this into \eqref{AppendixB5}, we obtain
\begin{equation}\label{AppendixB6}
\frac{d \tilde{\varphi} (\varepsilon)}{d \varepsilon}|_{\varepsilon =0} = -
\frac{\text{tr} \mathbf A}{3 \varphi_0}, \quad
\tilde{\varphi}(\varepsilon) = \varphi_0 - \frac{\text{tr} \mathbf A}{3 \varphi_0} \varepsilon + O(\varepsilon^2).
\end{equation}

Interestingly, this estimation is exact in the case of isotropic $\mathbf{A}$.
Moreover, a higher-order approximation of $\varphi$ is possible.
However, the higher-order approximation would be even less
accurate dealing with finite $\varepsilon$, which may occur in practice.

\section*{Appendix C: order of approximation}

Let us show that IFEBM is first order accurate.
Toward that end we consider a typical time interval $(t_n, t_{n+1})$; ${}^{n+1} \mathbf{C}$ and
${}^{n} \mathbf{C}_{\text{i}}$ are given and $\det {}^{n} \mathbf{C}_{\text{i}}  =1$.
Let $\mathbf{C}^{\text{exact}}_{\text{i}}$ be the exact solution to equation \eqref{puba44}
with the initial condition $\mathbf C_{\text{i}}|_{t=t_n} = {}^{n} \mathbf{C}_{\text{i}}$ and $\mathbf{C} \equiv {}^{n+1} \mathbf{C}$;
${}^{n+1} \mathbf{C}^{\text{IFEBM}}_{\text{i}}$ is the corresponding IFEBM-solution.
It is sufficient to show that for small $\Delta t$
\begin{equation}\label{AppendixC1}
{}^{n+1} \mathbf{C}^{\text{IFEBM}}_{\text{i}} -
\mathbf{C}^{\text{exact}}_{\text{i}}(t_{n+1}) = O((\Delta t)^2).
\end{equation}
We re-write \eqref{puba44} in the compact form
\begin{equation}\label{AppendixC2}
\dot{\mathbf{C}}_{\text{i}} (t)  = \mathbf{f} (\mathbf{C}_{\text{i}}(t)) \mathbf{C}_{\text{i}}, \quad
\text{tr}\mathbf{f} \equiv 0.
\end{equation}
Since $\mathbf{C}^{\text{exact}}_{\text{i}}(\cdot)$ is smooth, the mean value theorem implies that there exists $t^* \in (t_n,t_{n+1})$ such that
\begin{equation}\label{AppendixC3}
\mathbf{C}^{\text{exact}}_{\text{i}}(t_{n+1})  = {}^{n} \mathbf{C}_{\text{i}} + \Delta t \ \mathbf{f} (
\mathbf{C}^{\text{exact}}_{\text{i}}(t^*) ) \mathbf{C}^{\text{exact}}_{\text{i}}(t^*).
\end{equation}
Due to the smoothness of $\mathbf{f}$, we have
\begin{equation}\label{AppendixC4}
\mathbf{C}^{\text{exact}}_{\text{i}}(t_{n+1})  = {}^{n} \mathbf{C}_{\text{i}} + \Delta t \ \mathbf{f}
( {}^{n} \mathbf{C}_{\text{i}} ) \ {}^{n} \mathbf{C}_{\text{i}} + O((\Delta t)^2).
\end{equation}
Since $\det \mathbf{C}^{\text{exact}}_{\text{i}}(t_{n+1})  =1$ and
$\det {}^{n} \mathbf{C}_{\text{i}}  =1$, using the Jacobi formula we have
\begin{equation}\label{AppendixC5}
{}^{n} \mathbf{C}^{-1}_{\text{i}} : \big(\Delta t \
\mathbf{f}({}^{n} \mathbf{C}_{\text{i}}) {}^{n}
\mathbf{C}_{\text{i}} + O((\Delta t)^2)\big) = 0 .
\end{equation}
On the other hand, for the exact solution of the discretized equation \eqref{ModEBM} we have
\begin{equation}\label{AppendixC6}
 {}^{n+1} \mathbf{C}_{\text{i}}  = {}^{n} \mathbf{C}_{\text{i}} + \Delta t \
 \mathbf{f} ({}^{n+1} \mathbf{C}_{\text{i}} ) \ {}^{n+1} \mathbf{C}_{\text{i}} +
\tilde{\varphi} \ {}^{n+1} \mathbf{C}_{\text{i}}, \quad
 \det  {}^{n+1} \mathbf{C}_{\text{i}}  =1.
\end{equation}
Due to the smoothness of $\mathbf{f}$,
\begin{equation}\label{AppendixC7}
 {}^{n+1} \mathbf{C}_{\text{i}}  = {}^{n} \mathbf{C}_{\text{i}} + \Delta t \
 \mathbf{f} ({}^{n} \mathbf{C}_{\text{i}} ) \ {}^{n} \mathbf{C}_{\text{i}} +
  \tilde{\varphi} \ {}^{n} \mathbf{C}_{\text{i}} + O((\Delta t)^2).
\end{equation}
Similarly to \eqref{AppendixC5}, using the incompressibility condition $\det  {}^{n+1} \mathbf{C}_{\text{i}}  =1$, we obtain
\begin{equation}\label{AppendixC8}
{}^{n} \mathbf{C}^{-1}_{\text{i}} : \big(\Delta t \
\mathbf{f}({}^{n} \mathbf{C}_{\text{i}}) {}^{n}
\mathbf{C}_{\text{i}} +
\tilde{\varphi} \ {}^{n} \mathbf{C}_{\text{i}} + O((\Delta t)^2)\big) = 0 .
\end{equation}
Subtracting \eqref{AppendixC5} from \eqref{AppendixC8}, we conclude that $\tilde{\varphi} = O((\Delta t)^2)$.
Having this in mind and subtracting \eqref{AppendixC4} from \eqref{AppendixC7}, we obtain the desired estimation
\begin{equation}\label{AppendixC9}
{}^{n+1} \mathbf{C}_{\text{i}} - \mathbf{C}^{\text{exact}}_{\text{i}}(t_{n+1}) =  \tilde{\varphi} + O((\Delta t)^2) = O((\Delta t)^2).
\end{equation}
Note that the IFEBM solution differs slightly from ${}^{n+1} \mathbf{C}_{\text{i}}$,
since in EFEBM the parameter $\varphi$ is not identified exactly, but estimated
using equation \eqref{EstimOfz}. The error in the estimation of $\varphi$ is of order $O(\varepsilon^2) = O((\Delta t)^2)$.
Thus, the corresponding $\tilde{\varphi}$ is still $O((\Delta t)^2)$ and \eqref{AppendixC9} is still valid:
\begin{equation}\label{AppendixC10}
{}^{n+1} \mathbf{C}^*_{\text{i}} - \mathbf{C}^{\text{exact}}_{\text{i}}(t_{n+1}) = O((\Delta t)^2).
\end{equation}
Finally, the projection $\overline{(\cdot)}$
(cf. \eqref{UpdateCi_correct}) brings changes
which are only $O((\Delta t)^2)$: ${}^{n+1} \mathbf{C}^*_{\text{i}} -
\overline{{}^{n+1} \mathbf{C}^*_{\text{i}}} = O((\Delta t)^2)$
(cf. Appendix B in \cite{ShutovLandgraf}). Thus, IFEBM is first order accurate.

\section*{Appendix D: exact preservation of the w-invariance}

A general definition of the w-invariance under isochoric change of the reference configuration
is presented in \cite{WInvariance}.
In a simple formulation, the w-invariance of a material model says that for any
volume-preserving change of the reference configuration
there is a transformation of initial conditions ensuring that
the predicted Cauchy stresses are not affected by the reference change.
The w-invariance represents a generalized symmetry of the constitutive equations
indicating that the analyzed material exhibits fluid-like properties.

In the particular case of constitutive
equations \eqref{UseinAppD} and \eqref{puba4}, the w-invariance property
is formulated in the following way. Let $\mathbf{F}_0$ be arbitrary
second-rank tensor, such that $\det \mathbf{F}_0 =1$. Consider a prescribed history of the
right Cauchy-Green tensor $\mathbf{C}(t)$, $t \in [0,T]$ and the initial
condition $\mathbf C_{\text{i}}|_{t=0} = \mathbf C_{\text{i}}^0$.
Let $\mathbf C^{\text{new}}_{\text{i}}$ be a new solution of
\eqref{UseinAppD}, \eqref{puba4} corresponding to the new loading programm
$\mathbf C^{\text{new}} (t)$ and new initial conditions
\begin{equation}\label{AppendixD1}
\mathbf C^{\text{new}} (t) :=  \mathbf{F}^{-\text{T}}_0 \mathbf C (t) \mathbf{F}^{-1}_0, \quad
\mathbf C^{\text{new}}_{\text{i}}|_{t=0} =  \mathbf{F}^{-\text{T}}_0 \mathbf C_{\text{i}}^0 \mathbf{F}^{-1}_0.
\end{equation}
System \eqref{UseinAppD}, \eqref{puba4} is w-invariant if and only if the original and new solutions a related through
\begin{equation}\label{AppendixD2}
\mathbf C^{\text{new}}_{\text{i}} (t) =  \mathbf{F}^{-\text{T}}_0 \mathbf C_{\text{i}} (t) \mathbf{F}^{-1}_0.
\end{equation}

Let $^{n+1} \mathbf{C}$,
$^{n} \mathbf{C}_{\text{i}}$, and $\mathbf{F}_0$ be given. Denote by $^{n+1} \mathbf{C}_{\text{i}}$ and
$^{n+1} \mathbf{C}^{\text{new}}_{\text{i}}$ the IFEBM-solutions pertaining to the original and new inputs, respectively:
\begin{equation}\label{AppendixD3}
\big( ^{n+1} \mathbf{C}, ^{n} \mathbf{C}_{\text{i}} \big) \
\stackrel{\text{IFEBM}}{\mapsto}  \ {}^{n+1} \mathbf{C}_{\text{i}}, \quad
\big( \mathbf{F}^{-\text{T}}_0 \ {}^{n+1} \mathbf{C} \ \mathbf{F}^{-1}_0,
\mathbf{F}^{-\text{T}}_0 \ {}^{n} \mathbf{C}_{\text{i}} \ \mathbf{F}^{-1}_0 \big) \
\stackrel{\text{IFEBM}}{\mapsto}  \ ^{n+1} \mathbf{C}^{\text{new}}_{\text{i}}
\end{equation}
In analogy to the continuous case \eqref{AppendixD2}, the algorithm is said to preserve the w-invariance if
\begin{equation}\label{AppendixD4}
{}^{n+1}\mathbf{C}^{\text{new}}_{\text{i}} = \mathbf{F}^{-\text{T}}_0  \ {}^{n+1} \mathbf{C}^{\text{new}}_{\text{i}} \
\mathbf{F}^{-1}_0.
\end{equation}
Algorithms satisfying
\eqref{AppendixD4} are advantageous over algorithms which
violate this symmetry restriction (cf. the discussion in \cite{ShutovLandgraf}).
A straightforward (but tedious) way of proving \eqref{AppendixD4} is to substitute
the IFEBM equations into \eqref{AppendixD4} and to check the identity.
A more elegant proof is based on the observation that the IFEBM on the reference
configuration (cf. table \ref{tab0}) is equivalent to
the IFEBM on the current configuration (cf. table \ref{tab2}).
More precisely, let ${}^{n+1}\mathbf{F}$ be any second-rank tensor, such that
${}^{n+1}\mathbf{F}^{\text{T}} \ {}^{n+1}\mathbf{F} = {}^{n+1}\mathbf{C}$.
Recall that, according to \eqref{Eulerian3} and \eqref{Eulerian5corrected}, ${}^{n+1} {\mathbf B}_{\text{e}}$
is a unique function of the trial value ${}^{n+1} {\mathbf B}^{\text{trial}}_{\text{e}}$.
The following computation steps yield exactly the IFEBM-solution:
\begin{equation}\label{AppendixD5}
{}^{n+1} {\mathbf B}^{\text{trial}}_{\text{e}} =
{}^{n+1} {\mathbf F} \  {}^{n} \mathbf{C}^{-1}_{\text{i}} \ {}^{n+1} {\mathbf F}^{\text{T}}, \quad
{}^{n+1} {\mathbf B}^{\text{trial}}_{\text{e}}  \stackrel{\eqref{Eulerian3}, \eqref{Eulerian5}}{\mapsto}
{}^{n+1} {\mathbf B}_{\text{e}}, \quad
{}^{n+1} \mathbf{C}_{\text{i}} = {}^{n+1} {\mathbf F}^{\text{T}} \ {}^{n+1} {\mathbf B}^{-1}_{\text{e}}
\ {}^{n+1} {\mathbf F}.
\end{equation}
Now, for the new inputs ${}^{n+1} \mathbf{C}^{\text{new}} =
\mathbf{F}^{-\text{T}}_0 \ {}^{n+1} \mathbf{C} \ \mathbf{F}^{-1}_0$ and
${}^{n} \mathbf{C}^{\text{new}}_{\text{i}} =  \mathbf{F}^{-\text{T}}_0
\ {}^{n} \mathbf{C}_{\text{i}} \ \mathbf{F}^{-1}_0$ we may put
${}^{n+1} \mathbf{F}^{\text{new}} = {}^{n+1} \mathbf{F} \mathbf{F}^{-1}_0$.
Then
\begin{equation}\label{AppendixD6}
({}^{n+1} {\mathbf B}^{\text{trial}}_{\text{e}})^{\text{new}} =
{}^{n+1} {\mathbf F}^{\text{new}} \  ({}^{n} \mathbf{C}^{\text{new}}_{\text{i}})^{-1} \
({}^{n+1} {\mathbf F}^{\text{new}})^{\text{T}} = {}^{n+1} {\mathbf B}^{\text{trial}}_{\text{e}}.
\end{equation}
Since the trial values for the original and new inputs coincide, we have
\begin{equation}\label{AppendixD7}
{}^{n+1} {\mathbf B}^{\text{new}}_{\text{e}} =
{}^{n+1} {\mathbf B}_{\text{e}}, \quad
{}^{n+1} \mathbf{C}^{\text{new}}_{\text{i}} = ({}^{n+1} {\mathbf F}^{\text{new}})^{\text{T}} \
{}^{n+1} {\mathbf B}^{-1}_{\text{e}} \
{}^{n+1} {\mathbf F}^{\text{new}} =
\mathbf{F}^{-\text{T}}_0  \  {}^{n+1} \mathbf{C}_{\text{i}} \ \mathbf{F}^{-1}_0.
\end{equation}
That is exactly the required relation between ${}^{n+1}
\mathbf{C}^{\text{new}}_{\text{i}}$  and ${}^{n+1}\mathbf{C}_{\text{i}}$.
In the same way, the w-invariance can be proved for the 2IEBM.

\section*{Appendix E: MEBM and EM}

We consider the initial value problem
\begin{equation*}\label{difur}
\dot{\mathbf C_{\text{i}}} (t) = \mathbf{f} (\mathbf C_{\text{i}}(t), t) \mathbf C_{\text{i}}(t), \quad
\mathbf C_{\text{i}}(0)=\mathbf C_{\text{i}}^0, \quad \det(\mathbf
C_{\text{i}}^0)=1, \quad \text{where} \ \text{tr} \mathbf{f} \equiv 0.
\end{equation*}
Assume that ${}^n \mathbf C_{\text{i}}$ is known.
The classical Euler Backward method (EBM) is based on the equation
\begin{equation}\label{EBMClassic}
{}^{n+1} \mathbf C_{\text{i}}^{\text{EBM}} =  {}^n \mathbf C_{\text{i}} +
\Delta t \ \mathbf{f} ({}^{n+1}\mathbf C_{\text{i}}^{\text{EBM}}, t_{n+1}) \ {}^{n+1}\mathbf C_{\text{i}}^{\text{EBM}}.
\end{equation}
In the current study we use its modification, called modified Euler backward method (MEBM),
which guaranties that $\det(\mathbf {}^{n+1}\mathbf C_{\text{i}})=1$
\begin{equation}\label{MEBM345}
{}^{n+1} \mathbf C_{\text{i}}^{\text{MEBM}} = \overline{ {}^n \mathbf C_{\text{i}} + \Delta t
\mathbf{f} ({}^{n+1}\mathbf C_{\text{i}}^{\text{MEBM}}, t_{n+1}) \ {}^{n+1}\mathbf C_{\text{i}}^{\text{MEBM}} }.
\end{equation}
Another modification of the Euler backward method was presented in \cite{Helm2}, also
to enforce the inelastic incompressibility.
The exponential method (EM) corresponds to the equation
\begin{equation}\label{Expo}
{}^{n+1} \mathbf C_{\text{i}}^{\text{EM}} = \exp\big(\displaystyle \Delta t \
\mathbf{f} ({}^{n+1} \mathbf C_{\text{i}}^{\text{EM}}, t_{n+1})\big) \ {}^n \mathbf C_{\text{i}}.
\end{equation}
MEBM and EM exactly preserve the geometric property \eqref{geopro}:
${}^{n+1} \mathbf{C}_{\text{i}} \in Sym$, $\det({}^{n+1} \mathbf{C}_{\text{i}}) = 1$.

\end{document}